\documentclass[a4paper,headlines=2.1]{scrartcl}

\usepackage[utf8]{inputenc}
\usepackage{amsmath,amsthm,amssymb}
\usepackage{amsrefs}
\usepackage{tikz}
\usetikzlibrary{arrows,quotes,shapes,automata,petri,positioning,calc,arrows.meta,fit,chains}
\usepackage{ifthen}
\usetikzlibrary{decorations, decorations.markings} 
\usetikzlibrary{shapes.geometric}
\usetikzlibrary{tikzmark}
\usetikzlibrary{patterns,snakes}
\usepackage{hyperref}
\usepackage{bbm}
\usepackage{xcolor}
\usepackage{tabto}
\usepackage{bm}
\usepackage{enumitem}
\usepackage{authblk}
\usepackage{mathtools}

\usepackage{tabularx}
\newcolumntype{Y}{>{\centering\arraybackslash}X}

\usepackage{xr}

\allowdisplaybreaks

\newcommand{\N}{\mathbb{N}}
\newcommand{\E}{\mathbb{E}}

\newcommand{\cW}{\mathcal{W}}
\newcommand{\tZ}{\tilde{Z}}
\newcommand{\tA}{\tilde{A}}

\newcommand{\ol}{\overline}
\newcommand{\tr}{\operatorname{tr}}

\newcommand{\tb}{\tilde{b}}

\newcommand{\cS}{\mathcal{S}}
\renewcommand{\mod}{ \, \operatorname{mod} \, }

\tikzset{
	graydouble/.style={minimum height=0.4cm-\pgflinewidth, minimum width=2*0.4cm-\pgflinewidth, outer sep=0pt, fill=black!30},
	gray/.style={minimum height=0.4cm-\pgflinewidth, minimum width=0.4cm-\pgflinewidth, outer sep=0pt, fill=black!30},
	black/.style={minimum height=0.4cm-\pgflinewidth, minimum width=0.4cm-\pgflinewidth, outer sep=0pt, fill=black},
	white/.style={minimum size=0.4cm-\pgflinewidth, outer sep=0pt, fill=white, scale = \s, very thick},
	blue/.style={minimum size=0.4cm-\pgflinewidth, outer sep=0pt, fill=blue!20},
	red/.style={minimum size=0.4cm-\pgflinewidth, outer sep=0pt, fill=red!20},
}

\tikzset{
	-<-/.style args={#1 #2 #3}{
		decoration={
			markings,
			mark= at position 0.5 with
			{
				\ifthenelse{#3 = 1}
				{
					\fill[#2] (#1/6.0,0pt) -- (0.5*#1, #1/3.0) -- (-0.5*#1,0pt) -- (0.5*#1, #1/-3.0);   
				}
				{
					\ifthenelse{#3 = 2}
					{
						\fill[#2] (#1/-2.0,0pt) -- (0.5*#1, #1/3.0) -- (0.5*#1, #1/-3.0);   
					}
					{
						\ifthenelse{#3 = 3}
						{
							\draw[semithick, #2]  (0.533*#1,#1/2) -- (-0.433*#1, 0) -- (0.533*#1,-#1/2);  
						}{}
					}
				}
			},
		},
		postaction={decorate}
	},
	-<-/.default={6pt black 1}
}

\newsavebox\thesmashminipage

\definecolor{dgreen}{RGB}{0, 120, 0}

\newcommand\s{1}
\renewcommand\d{\s cm}
\tikzset{
	blacknode/.style={
		circle,
		thick,
		draw=black,
		fill=gray!50,
		minimum size=0.05*\d,
		scale = \s,
	},
	whitenode/.style={
		circle,
		thick,
		draw=black,
		minimum size=0.05*\d,
		scale = \s,
	},
	noborder/.style={
		circle,
		thick,
		minimum size=0.05*\d,
		scale = \s,
	},
	el/.style = {inner sep=2pt, align=left, sloped, font=\tiny},
	every label/.append style = {font=\tiny},
	position/.style args={#1:#2 from #3}{
		at=(#3.#1), anchor=#1+180, shift=(#1:#2)
	}
}

\numberwithin{equation}{section}

\AtBeginEnvironment{biblist}{\catcode`\#=12 }

\title{A recursion formula for mixed trace moments of isotropic Wishart matrices and the Gaussian unitary/orthogonal ensembles}
\date{}
\author{Ben Deitmar}
\affil{\small \textit{Department of Mathematical Stochastics, ALU Freiburg \protect\\ Ernst-Zermelo-Str. 1, 79104 Freiburg, Germany \protect\\ E-mail: ben.deitmar@stochastik.uni-freiburg.de}}
\begin{document}
	\thispagestyle{empty}
	\maketitle
	\vspace{-1.2cm}
	
	\begin{abstract}
		\begin{center}\textbf{Abstract}\end{center}
		\noindent
		Exact recursion formulas for mixed moments of four fundamental random matrix ensembles are derived. The reason such recursive formulas are possible is closely related to properties of polygon gluings studied by Harer and Zagier as well as Akhmedov and Shakirov. The proofs of the formulas are direct and written in such a way that they do not rely on understanding of polygon gluings.
	\end{abstract}
	
	

	\section{Introduction}
	Let $\bm{X}$ be a $(p \times n)$ random matrix with iid standard normal entries and let $\bm{Y}$ be a $(p \times n)$ random matrix with iid standard complex normal entries. Further let $\bm{Z}$ be an $(n \times n)$ random hermitian matrix, where the entries $(Z_{i,j})_{i \leq j}$ are independent and $Z_{i,i} \sim \mathcal{N}(0,1)$ for all $i \leq n$ as well as $Z_{i,j} = \ol{Z}_{j,i} \sim \mathcal{CN}(0,1)$ for all $i<j \leq n$. Similarly, let $\tilde{\bm{Z}}$ be an $(n \times n)$ random symmetric matrix, where the entries $(\tZ_{i,j})_{i \leq j}$ are independent and $\tZ_{i,i} \sim \mathcal{N}(0,2)$ for all $i \leq n$ as well as $\tZ_{i,j} = \tZ_{j,i} \sim \mathcal{N}(0,1)$ for all $i<j \leq n$.\\
	\\
	The matrices $\frac{1}{\sqrt{n}}\bm{Z}$ and $\frac{1}{\sqrt{n}} \tilde{\bm{Z}}$ are of the Gaussian unitary and Gaussian orthogonal ensemble respectively. The matrices $\frac{1}{n} \bm{X} \bm{X}^T$ and $\frac{1}{n} \bm{Y} \bm{Y}^*$ are real and complex isotropic Wishart matrices.\\
	\\
	In this paper we give exact recursive formulas for mixed trace moments of the type
	\begin{align*}
		& E_X(l) := \E\bigg[ \prod\limits_{k=1}^K \tr\big( (\bm{X}\bm{X}^T)^{l_k} \big) \bigg] \ \ , \ \ E_Y(l) := \E\bigg[ \prod\limits_{k=1}^K \tr\big( (\bm{Y}\bm{Y}^*)^{l_k} \big) \bigg]
	\end{align*}
	\hspace{-0.7cm}\rule[-0.5cm]{8cm}{0.2mm}\\
	\\
	\hspace{1cm}
	Supported by the DFG Research Unit 5381
	
	\newpage
	\noindent
	and
	\begin{align*}
		& E_Z(l) := \E\bigg[ \prod\limits_{k=1}^K \tr\big( \bm{Z}^{l_k} \big) \bigg] \ \ , \ \ E_{\tZ}(l) := \E\bigg[ \prod\limits_{k=1}^K \tr\big( \tilde{\bm{Z}}^{l_k} \big) \bigg]
	\end{align*}
	for arbitrary $p,n,K \in \N$ and \textit{layouts} $l = (l_1,...,l_K) \in \N_0^K$.\\
	\\
	Assuming $l_1 > 0$, the formulas are as follows.
	\begin{itemize}
		\item[1)] Gaussian unitary ensemble:
		\begin{align*}
			& \E\bigg[ \prod\limits_{r=1}^K \tr\big( \bm{Z}^{l_r} \big) \bigg] = \sum\limits_{q=1}^{l_1-1} \E\bigg[ \tr\big( \bm{Z}^{q-1} \big) \tr\big( \bm{Z}^{l_1-q-1} \big) \prod\limits_{r=2}^K \tr\big( \bm{Z}^{l_r} \big) \bigg]\\
			& \hspace{3cm} + \sum\limits_{k=2}^K l_k \E\bigg[ \tr\big( \bm{Z}^{l_1+l_k-2} \big) \prod\limits_{\substack{r=1 \\ r \neq k}}^K \tr\big( \bm{Z}^{l_r} \big) \bigg]
		\end{align*}
		For a complete formulation and proof see Theorem \ref{Thm_GUE_Recursion}.
		
		\item[2)] Gaussian orthogonal ensemble:
		\begin{align*}
			& \E\bigg[ \prod\limits_{r=1}^K \tr\big( \tilde{\bm{Z}}^{l_r} \big) \bigg] = (l_1-1) \E\bigg[ \tr\big( \tilde{\bm{Z}}^{l_1-2} \big) \prod\limits_{r=2}^K \tr\big( \tilde{\bm{Z}}^{l_r} \big) \bigg]\\
			& \hspace{3cm} + \sum\limits_{q=1}^{l_1-1} \E\bigg[ \tr\big( \tilde{\bm{Z}}^{q-1} \big) \tr\big( \tilde{\bm{Z}}^{l_1-q-1} \big) \prod\limits_{r=2}^K \tr\big( \tilde{\bm{Z}}^{l_r} \big) \bigg]\\
			& \hspace{3cm} + 2\sum\limits_{k=2}^K l_k \E\bigg[ \tr\big( \tilde{\bm{Z}}^{l_1+l_k-2} \big) \prod\limits_{\substack{r=1 \\ r \neq k}}^K \tr\big( \tilde{\bm{Z}}^{l_r} \big) \bigg]
		\end{align*}
		For a complete formulation and proof see Theorem \ref{Thm_GOE_Recursion}.

		\item[3)] Complex isotropic Wishart ensemble:
		\begin{align*}
			& \E\bigg[ \prod\limits_{r=1}^K \tr\big( (\bm{Y}\bm{Y}^*)^{l_r} \big) \bigg] = n \E\bigg[ \tr\big( (\bm{Y}\bm{Y}^*)^{l_1-1} \big) \prod\limits_{r=2}^K \tr\big( (\bm{Y}\bm{Y}^*)^{l_r} \big) \bigg]\\
			& \hspace{4cm} + \sum\limits_{q=1}^{l_1-1} \E\bigg[ \tr\big( (\bm{Y}\bm{Y}^*)^{q} \big) \tr\big( (\bm{Y}\bm{Y}^*)^{l_1-q-1} \big) \prod\limits_{r=2}^K \tr\big( (\bm{Y}\bm{Y}^*)^{l_r} \big) \bigg]\\
			& \hspace{4cm} + \sum\limits_{k=2}^K l_k \E\bigg[ \tr\big( (\bm{Y}\bm{Y}^*)^{l_1+l_k-1} \big) \prod\limits_{\substack{r=1 \\ r \neq k}}^K \tr\big( (\bm{Y}\bm{Y}^*)^{l_r} \big) \bigg]
		\end{align*}
		For a complete formulation and proof see Theorem \ref{Thm_WC_Recursion}.

		\item[4)] Real isotropic Wishart ensemble:
		\begin{align*}
			& \E\bigg[ \prod\limits_{r=1}^K \tr\big( (\bm{X}\bm{X}^T)^{l_r} \big) \bigg] = (n+l_1-1) \E\bigg[ \tr\big( (\bm{X}\bm{X}^T)^{l_1-1} \big) \prod\limits_{r=2}^K \tr\big( (\bm{X}\bm{X}^T)^{l_r} \big) \bigg]\\
			& \hspace{4cm} + \sum\limits_{q=1}^{l_1-1} \E\bigg[ \tr\big( (\bm{X}\bm{X}^T)^{q} \big) \tr\big( (\bm{X}\bm{X}^T)^{l_1-q-1} \big) \prod\limits_{r=2}^K \tr\big( (\bm{X}\bm{X}^T)^{l_r} \big) \bigg]\\
			& \hspace{4cm} + 2\sum\limits_{k=2}^K l_k \E\bigg[ \tr\big( (\bm{X}\bm{X}^T)^{l_1+l_k-1} \big) \prod\limits_{\substack{r=1 \\ r \neq k}}^K \tr\big( (\bm{X}\bm{X}^T)^{l_r} \big) \bigg]
		\end{align*}
		For a complete formulation and proof see Theorem \ref{Thm_WR_Recursion}.
	\end{itemize}
	The recursive formulas are proved in order of increasing difficulty, so Theorem \ref{Thm_GUE_Recursion} proves (1), Theorem \ref{Thm_WC_Recursion} proves (3), Theorem \ref{Thm_WR_Recursion} proves (4) and Theorem \ref{Thm_GOE_Recursion} proves (2).\\
	\\
	The only random matrix ensembles for which explicit and non-recursive formulas for mixed moments are known are the Haar unitary and the Haar orthogonal ensembles. Diaconis and Evans found the explicit formulas in \cite{DiaconisEvans} under an assumption roughly equivalent to $n \geq l_1+...+l_K$. Achieving non-recursive formulas for the ensembles studied here would prove extremely useful in the study of their spectral properties. For example the asymptotic behavior of
	\begin{align*}
		& \E\big[ \tr\big( \bm{Z}^{\lfloor t_1 n^{\frac{2}{3}}\rfloor} \big) \cdots \tr\big( \bm{Z}^{\lfloor t_K n^{\frac{2}{3}}\rfloor} \big) \big]
	\end{align*}
	is known to determine the behavior of the largest eigenvalue of $\bm{Z}$ at the Tracy-Widom scale, which is why such mixed trace moments were already extensively studied in \cite{Peche}, \cite{Soshnikov1}, \cite{Soshnikov2}, \cite{Soshnikov3} and \cite{Soshnikov4}. The strong universality results therein imply that a non-recursive result in any of the four cases studied here would already yield a greater understanding of Tracy-Widom laws in much more general settings. The fact that Akhmedov and Shakirov were able to find a non-recursive formula to a similar, but still different recursion in \cite{AkhmedovShakirov} suggests that non-recursive solutions may be achievable.\\
	\\
	Explicit formulas for the (non-mixed) singular trace moments $\E\big[ \tr\big( \bm{Z}^{2m} \big) \big]$ as well as $\E\big[ \tr\big( (\bm{X}\bm{X}^T)^{m} \big) \big]$ and $\E\big[ \tr\big( (\bm{Y}\bm{Y}^*)^{m} \big) \big]$ are known thanks to Haarer-Zagier in \cite{HarerZagier} and E. Vassilieva in \cite{Vassilieva}. The formula in the GUE case is not written explicitly in \cite{HarerZagier}, but following their work we give the formula in (\ref{Eq_HarerZagierExplicit}). The formulas for the real and complex Wishart case can be found in Corollaries 1.8 and 1.9 of \cite{Vassilieva}. For the singular trace moments in the GOE case $\E\big[ \tr\big( \tilde{\bm{Z}}^{2m} \big) \big]$ a five-term recurrence relation is known thanks to M. Ledoux in \cite{LedouxRecursion}.
	\\
	\\
	The recursion from Theorem \ref{Thm_GUE_Recursion} for the GUE is not new, as it also follows from Section 5 in \cite{HarerZagier} by Harer and Zagier. We however still prove it as a precursor to Theorems \ref{Thm_WC_Recursion}, \ref{Thm_WR_Recursion} and \ref{Thm_GOE_Recursion}. The real Wishart case from Theorem \ref{Thm_WR_Recursion} was already studied by Pielaszkiewicz, Von Rosen and Singull in \cite{WishartRecursion}, though they made a minor error in their formula. To the best of our knowledge the results of Theorems \ref{Thm_WC_Recursion} and \ref{Thm_GOE_Recursion} for the complex Wishart ensemble and the GOE respectively are entirely new.

	\subsection{Trace moments and polygon gluings}\label{Subsection_Polygongluings}
	There is a well-known link between trace moments of the GUE and the number of ways to glue to edges of a polygon together to receive an orientable surface of a given genus. This link was first described by Harer and Zagier in \cite{HarerZagier} to prove the recursive formula
	\begin{align}\label{Eq_HarerZagierRecusrion}
		& (m+1) \E\big[ \tr(\bm{Z}^{2m}) \big] \nonumber\\
		& = (4m-2)n \E\big[ \tr(\bm{Z}^{2m-2}) \big] + (m-1)(2m-1)(2m-3) \E\big[ \tr(\bm{Z}^{2m-4}) \big] \ ,
	\end{align}
	which allows for the explicit representation
	\begin{align}\label{Eq_HarerZagierExplicit}
		& \E\big[ \tr(\bm{Z}^{2m}) \big] = \frac{(2m)!}{m!} \sum\limits_{r=1}^{n \land (m+1)} \frac{{n \choose r} {m \choose r-1}}{2^{m+1-r}} \ .
	\end{align}
	These equalities were first found by Harer and Zagier and later simplified by Haagerup and Thorbj\o rnsen in \cite{Haagerup}, but the more readable formulation for (\ref{Eq_HarerZagierRecusrion}) is from Theorem 1 in \cite{LedouxRecursion} by M. Ledoux, where he goes on to find an analogous recursion formula for $\E\big[ \tr(\tilde{\bm{Z}}^{2m}) \big]$ in the GOE case.\\
	\\
	Akhmedov and Shakirov took this approach further in \cite{AkhmedovShakirov} by considering 'incomplete' gulings of a polygon such that the unglued edges form polygons of given lengths $(l_1,...,l_K)$. They were able to give a recursive formula for the number of such incomplete gluings with a given genus (for a certain definition of genus). Unfortunately, although they were able to derive the explicit expression from their recursive formula, their results do not seem applicable to random matrix theory and their solution to a similar recursion does not seem applicable for solving the recursive formulas derived in this article.\\
	\\
	In order to find formulas for mixed trace moments of the GUE for a layout $l=(l_1,...,l_K)$, one instead would need to count the number of ways to glue the edges of multiple polygons $P_{l_1},...,P_{l_K}$ such that the resulting oriented surface is of a given genus. Here a possible pairing $\pi = \big\{ \{e^1_1,e^1_5\},\{e^1_2,e^2_1\},... \big\}$ of edges for the gluing of $P_{8},P_{5},...,P_{6}$:
	\begin{align*}
		\renewcommand{\s}{1}
		\renewcommand{\d}{1cm}
		& \begin{array}{l} 
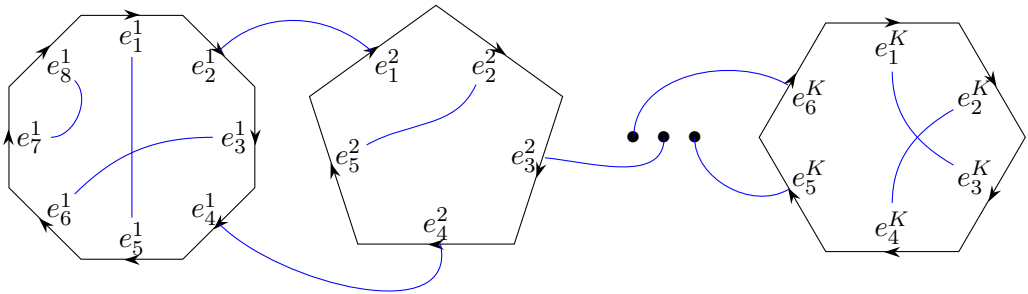
\begin{tikzpicture}
				\node (pol1) [draw, minimum size=3.5*\d, regular polygon, regular polygon sides=8,
				]  at (-4,0) {};
				\foreach \x/\y/\i in {1/2/1,2/3/{8},3/4/7,4/5/6,5/6/5,6/7/4,7/8/3,8/1/2}
				\path[auto=left, -<-]
				(pol1.corner \x)--(pol1.corner \y)
				node[midway, circle, inner sep=0cm, draw=none](e1_\i){$e^1_ {\i}$};
				\node (pol2) [draw, minimum size=3.5*\d, regular polygon, regular polygon sides=5] at (0,0) {};
				\foreach \x/\y/\i in {1/2/1,2/3/{5},3/4/4,4/5/3,5/1/2}
				\path[auto=left, -<-]
				(pol2.corner \x)--(pol2.corner \y)
				node[midway, circle, inner sep=0cm, draw=none](e2_\i){$e^2_{\i}$};
				\node at (3,0) {$\bullet \ \bullet \ \bullet$};
				\node (polK) [draw, minimum size=3.5*\d, regular polygon, regular polygon sides=6] at (6,0) {};
				\foreach \x/\y/\i in {1/2/1,2/3/{6},3/4/5,4/5/4,5/6/3,6/1/2}
				\path[auto=left, -<-]
				(polK.corner \x)--(polK.corner \y)
				node[midway, circle, inner sep=0cm, draw=none](eK_\i){$e^K_ {\i}$};
				\path[-]
				(e1_1) edge [bend right=0, color=blue]  (e1_5)
				(e1_7) edge [out=-45-45*7,in=-45-45*8, color=blue] (e1_8)
				(e1_6) edge [out=-45-45*6,in=-45-45*3, color=blue] (e1_3)
				(e1_2) edge [out=-45-45*2+180,in=30-72*1+180, color=blue] (e2_1)
				(e1_4) edge [out=-45-45*4+180,in=30-72*4+180, color=blue] (e2_4)
				(e2_5) edge [out=30-72*5,in=30-72*2, color=blue] (e2_2)
				(e2_3) edge [out=30-72*3+180,in=-90, color=blue] (3,0)
				(2.6,0) edge [out=90,in=-30-60*6+180, color=blue] (eK_6)
				(3.4,0) edge [out=-90,in=-30-60*5+180, color=blue] (eK_5)
				(eK_1) edge [out=-30-60*1,in=-30-60*3, color=blue] (eK_3)
				(eK_2) edge [out=-30-60*2,in=-30-60*4, color=blue] (eK_4)
				;
		\end{tikzpicture} \end{array}
	\end{align*}
	Euler's characteristic formula guarantees
	\begin{align*}
		& 2-2g = V - L/2 + K \ ,
	\end{align*}
	where $g$ is the genus of the resultant surface and $V$ is the number of vertices remaining in the ribbon graph determined by the pairing. Define $V_{L,K}(g) = 2-2g+\frac{L}{2}-K$ and $g_{L,K}(V) = \frac{L}{4}+1-\frac{V+K}{2}$.\\
	\\
	Let $\varepsilon_{g}(l_1,...,l_K)$ denote the number of different gluings/pairings of polygons $P_{l_1},...,P_{l_K}$ such that the resulting orientable surface is of genus $g$, then it is provable by Wick's Theorem (aka. Isserlis' Theorem) that
	\begin{align}\label{Eq_GUEOrderExpansion}
		& E_Z(l) = \sum\limits_{g=0}^{\lfloor g_{L,K}(0) \rfloor} n^{V_{L,K}(g)} \varepsilon_{g}(l_1,...,l_K)
	\end{align}
	under the assumption $l_1,...,l_K > 0$.\\
	\\
	A recursion formula for $E_Z(l)$ can then be found by mapping such polygon gluings onto gluings of polygons with two less total edges. If the edge $e^k_q$ paired with $e^1_1$ is still part of the same polygon, i.e. $k=1$, then we can contract along the connection of the edges before removing both such that we are left with a gluing of polygons $P_{q-2},P_{l_1-q},P_{l_2},...,P_{l_k}$. Here a picture of what the above pairing would look like after removing $e^1_1$ and its paired edge $e^1_5$:
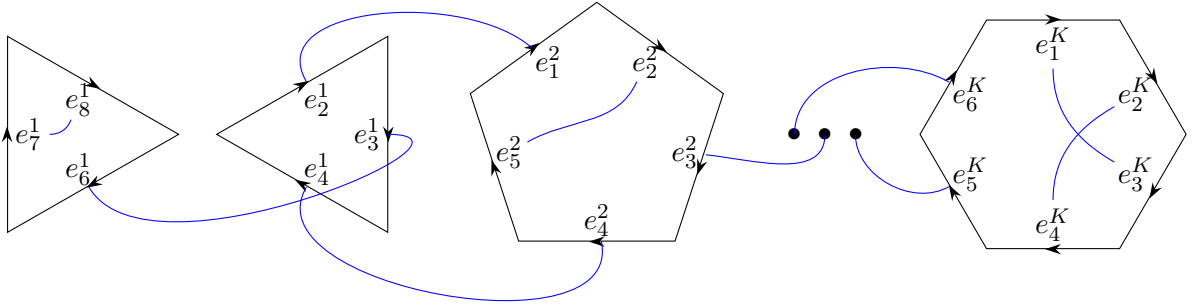
	\begin{align*}
		\renewcommand{\s}{1}
		\renewcommand{\d}{1cm}
		& \hspace{-0.5cm} \begin{array}{l} \begin{tikzpicture}
				\node (pol11) [draw, minimum size=3*\d, regular polygon, regular polygon sides=3,
				rotate=30]  at (-7,0) {};
				\foreach \x/\y/\i in {1/2/7,2/3/6,3/1/8}
				\path[auto=left, -<-]
				(pol11.corner \x)--(pol11.corner \y)
				node[midway, circle, inner sep=0cm, draw=none](e1_\i){$e^1_ {\i}$};
				\node (pol12) [draw, minimum size=3*\d, regular polygon, regular polygon sides=3,
				rotate=-30]  at (-3.5,0) {};
				\foreach \x/\y/\i in {1/2/2,2/3/4,3/1/3}
				\path[auto=left, -<-]
				(pol12.corner \x)--(pol12.corner \y)
				node[midway, circle, inner sep=0cm, draw=none](e1_\i){$e^1_ {\i}$};
				\node (pol2) [draw, minimum size=3.5*\d, regular polygon, regular polygon sides=5] at (0,0) {};
				\foreach \x/\y/\i in {1/2/1,2/3/{5},3/4/4,4/5/3,5/1/2}
				\path[auto=left, -<-]
				(pol2.corner \x)--(pol2.corner \y)
				node[midway, circle, inner sep=0cm, draw=none](e2_\i){$e^2_{\i}$};
				\node at (3,0) {$\bullet \ \bullet \ \bullet$};
				\node (polK) [draw, minimum size=3.5*\d, regular polygon, regular polygon sides=6] at (6,0) {};
				\foreach \x/\y/\i in {1/2/1,2/3/{6},3/4/5,4/5/4,5/6/3,6/1/2}
				\path[auto=left, -<-]
				(polK.corner \x)--(polK.corner \y)
				node[midway, circle, inner sep=0cm, draw=none](eK_\i){$e^K_ {\i}$};
				\path[-]
				(e1_7) edge [out=-45-45*7,in=-110, color=blue] (e1_8)
				(e1_6) edge [out=-60,in=0, color=blue] (e1_3)
				(e1_2) edge [out=120,in=30-72*1+180, color=blue] (e2_1)
				(e1_4) edge [out=-120,in=30-72*4+180, color=blue] (e2_4)
				(e2_5) edge [out=30-72*5,in=30-72*2, color=blue] (e2_2)
				(e2_3) edge [out=30-72*3+180,in=-90, color=blue] (3,0)
				(2.6,0) edge [out=90,in=-30-60*6+180, color=blue] (eK_6)
				(3.4,0) edge [out=-90,in=-30-60*5+180, color=blue] (eK_5)
				(eK_1) edge [out=-30-60*1,in=-30-60*3, color=blue] (eK_3)
				(eK_2) edge [out=-30-60*2,in=-30-60*4, color=blue] (eK_4)
				;
		\end{tikzpicture} \end{array}
	\end{align*}
	The edges would of course have to be renamed and we have swapped the positions of the first two polygons such that the effects of contracting along the connection $\{e^1_1,e^1_5\}$ are more obvious.\\
	\\
	Similarly, when the edge $e^k_q$ paired with $e^1_1$ is not in the first polygon, i.e. $k>1$, contraction along the pairing and removal of the two edges creates one new polygon $P_{l_1+l_k-2}$ from the two polygons $P_{l_1}$ and $P_{l_k}$. The rest of the pairing remains a valid gluing for the polygons $P_{l_1+l_k-2},P_2,...,\widehat{P}_{l_k},...,P_{l_K}$. Here an example before the contraction and removal:
	\begin{align*}
		\renewcommand{\s}{1}
		\renewcommand{\d}{1cm}
		& \begin{array}{l} 
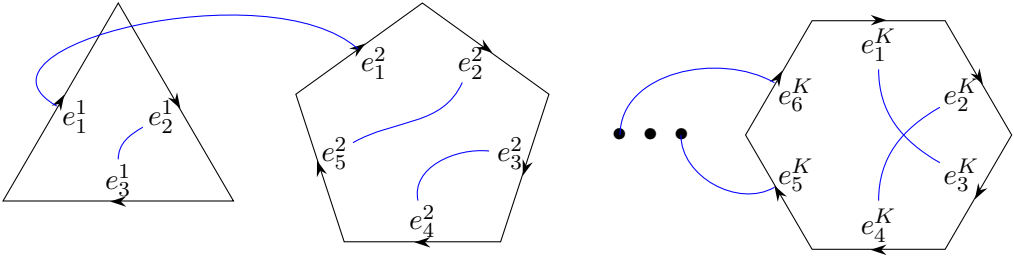
\begin{tikzpicture}
				\node (pol1) [draw, minimum size=3.5*\d, regular polygon, regular polygon sides=3,
				]  at (-4,0) {};
				\foreach \x/\y/\i in {1/2/1,2/3/3,3/1/2}
				\path[auto=left, -<-]
				(pol1.corner \x)--(pol1.corner \y)
				node[midway, circle, inner sep=0cm, draw=none](e1_\i){$e^1_ {\i}$};
				\node (pol2) [draw, minimum size=3.5*\d, regular polygon, regular polygon sides=5] at (0,0) {};
				\foreach \x/\y/\i in {1/2/1,2/3/{5},3/4/4,4/5/3,5/1/2}
				\path[auto=left, -<-]
				(pol2.corner \x)--(pol2.corner \y)
				node[midway, circle, inner sep=0cm, draw=none](e2_\i){$e^2_{\i}$};
				\node at (3,0) {$\bullet \ \bullet \ \bullet$};
				\node (polK) [draw, minimum size=3.5*\d, regular polygon, regular polygon sides=6] at (6,0) {};
				\foreach \x/\y/\i in {1/2/1,2/3/{6},3/4/5,4/5/4,5/6/3,6/1/2}
				\path[auto=left, -<-]
				(polK.corner \x)--(polK.corner \y)
				node[midway, circle, inner sep=0cm, draw=none](eK_\i){$e^K_ {\i}$};
				\path[-]
				(e1_1) edge [out=90-120*1+180,in=30-72*1+180, color=blue]  (e2_1)
				(e1_3) edge [out=90-120*3,in=90-120*2, color=blue] (e1_2)
				(e2_3) edge [out=30-72*3,in=30-72*4, color=blue] (e2_4)
				(e2_5) edge [out=30-72*5,in=30-72*2, color=blue] (e2_2)
				(2.6,0) edge [out=90,in=-30-60*6+180, color=blue] (eK_6)
				(3.4,0) edge [out=-90,in=-30-60*5+180, color=blue] (eK_5)
				(eK_1) edge [out=-30-60*1,in=-30-60*3, color=blue] (eK_3)
				(eK_2) edge [out=-30-60*2,in=-30-60*4, color=blue] (eK_4)
				;
		\end{tikzpicture} \end{array}
	\end{align*}
	and after the contraction and removal:
	\begin{align*}
		\renewcommand{\s}{1}
		\renewcommand{\d}{1cm}
		& \begin{array}{l} 
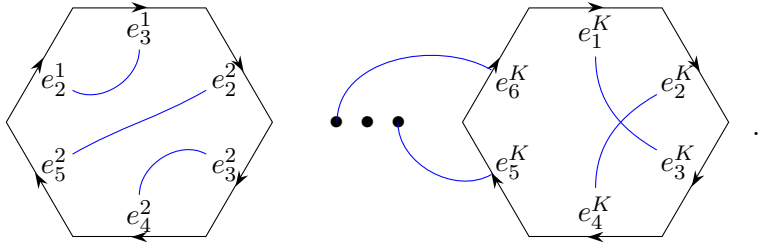
\begin{tikzpicture}
				\node (pol2) [draw, minimum size=3.5*\d, regular polygon, regular polygon sides=6] at (0,0) {};
				\foreach \x/\y/\c/\i in {1/2/1/3,2/3/1/2,3/4/2/5,4/5/2/4,5/6/2/3,6/1/2/2}
				\path[auto=left, -<-]
				(pol2.corner \x)--(pol2.corner \y)
				node[midway, circle, inner sep=0cm, draw=none](e\c_\i){$e^{\c}_{\i}$};
				\node at (3,0) {$\bullet \ \bullet \ \bullet$};
				\node (polK) [draw, minimum size=3.5*\d, regular polygon, regular polygon sides=6] at (6,0) {};
				\foreach \x/\y/\i in {1/2/1,2/3/{6},3/4/5,4/5/4,5/6/3,6/1/2}
				\path[auto=left, -<-]
				(polK.corner \x)--(polK.corner \y)
				node[midway, circle, inner sep=0cm, draw=none](eK_\i){$e^K_ {\i}$};
				\path[-]
				(e1_3) edge [out=90-60*3,in=90-60*2, color=blue] (e1_2)
				(e2_3) edge [out=90-60*5,in=90-60*6, color=blue] (e2_4)
				(e2_5) edge [out=90-60*1,in=90-60*4, color=blue] (e2_2)
				(2.6,0) edge [out=90,in=-30-60*6+180, color=blue] (eK_6)
				(3.4,0) edge [out=-90,in=-30-60*5+180, color=blue] (eK_5)
				(eK_1) edge [out=-30-60*1,in=-30-60*3, color=blue] (eK_3)
				(eK_2) edge [out=-30-60*2,in=-30-60*4, color=blue] (eK_4)
				;
		\end{tikzpicture} \end{array} \ .
	\end{align*}
	This method of removing edges does not remove vertices of the ribbon graph, while it can happen that some of the new polygons are empty (i.e. they have length $l_k=0$). This only becomes a problem, if the first polygon is empty, since then there is no edge $e^1_1$ to use for the next recursive step of contraction and removal. We thus discard the first polygon, if it is empty. As this also constitutes a removal of an isolated vertex from the ribbon graph, we must keep a tally of how many vertices were removed this way. Luckily, in the recursive formula this tally is simply realized with an additional factor $n$.\\
	\\
	We have thus far heuristically argued the plausibility of Theorem \ref{Thm_GUE_Recursion}.

	\subsection{Definition (Multigraphs by route)}\label{Def_Graphs}
	For a vertex set $V = [n] = \{1,...,n\}$ and a give \textit{route} $\bm{i}^1 \in [n]^m$ we define the directed multi-graph $G_{\bm{i}^1} = (V,E,f)$, with edge set $E = (e^1_1,...,e^1_m)$ and
	\begin{align*}
		& f: E \rightarrow V \times V \ \ ; \ \ e^1_q \mapsto (i^1_q,i^1_{(q \mod m)+1}) \ .
	\end{align*}
	By construction the route $\bm{i}^1$ is an Euler path on $G_{\bm{i}^1}$. Let $A(G_{\bm{i}^1}) \in \N_0^{n \times n}$ denote the adjacency matrix of $G_{\bm{i}^1}$, that is
	\begin{align*}
		& A_{v,w}(G_{\bm{i}^1}) = \#\{e^1_q \in E \mid f(e^1_q) = (v,w)\} \ .
	\end{align*}
	We can also define the undirected multi-graph $\tilde{G}_{\bm{i}^1} = (V,E,\tilde{f})$ to a given route $\bm{i}^1 \in [n]^m$ with the same edges, but ignorant of the edges direction. In this case we can define $\tilde{f}$ by
	\begin{align*}
		& \tilde{f} : E \rightarrow \{M \subset V \mid \#M \leq 2\} \ \ ; \ \ e^1_q \mapsto \{i^1_q,i^1_{(q \mod m)+1}\} \ .
	\end{align*}
	(If $\#f(e^1_q) = 1$, the edge $e^1_q$ is a self-loop.) Again $A(\tilde{G}_{\bm{i}^1}) \in \N_0^{n \times n}$ denotes the adjacency matrix of $\tilde{G}_{\bm{i}^1}$, i.e. 
	\begin{align*}
		& A_{v,w}(\tilde{G}_{\bm{i}^1}) = \#\{e^1_q \in E \mid f(e^1_q) = \{v,w\}\} \ . 
	\end{align*}
	Note that for this definition of the adjacency matrix self-loops are only counted once.
	For a sequence of routes $\bm{i} = (\bm{i}^1,...,\bm{i}^K) \in \bigtimes\limits_{k=1}^K [n]^{l_k}$ define the joined adjacency matrices
	\begin{align}\label{Eq_DefJoinedAdjacency}
		& A(\bm{i}) := \sum\limits_{k=1}^K A(G_{\bm{i}^k}) \ \ ; \ \ \tilde{A}(\bm{i}) := \sum\limits_{k=1}^K A(\tilde{G}_{\bm{i}^k}) \ .
	\end{align}
	In Section \ref{SectionWishart} the vertex set $V$ will be $[p+n]$, so we must exchange all above occurrences of $n$ with $p+n$.

	\section{Gaussian Unitary Ensemble}
	
	\subsection{Theorem (Recursion for mixed trace moments of the GUE)}\label{Thm_GUE_Recursion}
	For $K \in \N$ and any $l=(l_1,...,l_{K}) \in \N_0^{K}$ - with even $L := l_1+...+l_K$ - the following recursive properties hold.
	\begin{itemize}
		\item[a)]
		If $l_1=0$, then $E_{Z}(l) = nE_{Z}(S_0(l))$ for
		\begin{align}\label{EqGUE_DefS0}
			& S_0(0,l_2,...,l_K) := (l_2,...,l_K) \ .
		\end{align}
		
		\item[b)]
		If $l_1>0$, then
		\begin{align*}
			& E_Z(l) = \sum\limits_{q=1}^{l_1-1} E_{Z}(S_{2,q}(l)) + \sum\limits_{k=2}^{K} l_k E_{Z}(S_{3,k}(l)) \ ,
		\end{align*}
		where
		\begin{align}
			& S_{2,q}(l) := \big(q-1,l_1-q-1,l_2,...,l_{K}\big) \label{EqGUE_DefS2}\\
			& S_{3,k}(l) := \big(l_1+l_k-2,l_2,...,\widehat{l_k},...,l_{K}\big) \ . \label{EqGUE_DefS3}
		\end{align}
	\end{itemize}
	Since (a) reduces $K$ by one and (b) reduces $L = l_1+...+l_K$ by two, the starting point
	\begin{align*}
		& E_{Z}\big( ((,)) \big) = \E[\tr(\bm{Z})^0] = n
	\end{align*}
	for $(K,L) = (1,0)$ allows for the calculation of all mixed moments, when $L$ is even. When $L$ is odd, it is easy to see that the mixed trace moment must be zero.
	\begin{proof}\
		\\
		The property (a) holds trivially, since the empty matrix product is here the $(n \times n)$ identity matrix. It remains to prove property (b).\\
		\\
		For now assume $l_1,...,l_{K}>0$, then expanding the sums in $\tr(\bm{Z}^{l_i})$ yields
		\begin{align*}
			& E_{Z}(l) = \sum\limits_{\substack{\bm{i}^1 \in [n]^{l_1} \\ : \\ \bm{i}^{K} \in [n]^{K}}} \E\bigg[ \prod\limits_{k=1}^{K} Z_{i^k_1,i^k_2}  \cdots Z_{i^k_{l_k-1},i^k_{l_k}} \, Z_{i^k_{l_k},i^k_1} \bigg] \ .
		\end{align*}
		Extracting the first entries of each $\bm{i}^k$ and summing over them first gives
		\begin{align*}
			& E_{Z}(l) = \sum\limits_{b_1,...,b_K \in [n]} \sum\limits_{\substack{\bm{i} \in \bigtimes\limits_{k=1}^K [n]^{l_k} \\ \forall k \leq K: \, i^k_1 = b_k}} \E\bigg[ \prod\limits_{k=1}^{K} Z_{i^k_1,i^k_2}  \cdots Z_{i^k_{l_k-1},i^k_{l_k}} \, Z_{i^k_{l_k},i^k_1} \bigg] \ .
		\end{align*}
		For easier notation define
		\begin{align*}
			& \mathbb{I}_{l,b} := \big\{ \bm{i} = (\bm{i}^1,...,\bm{i}^K) \in \bigtimes\limits_{k=1}^K [n]^{l_k} \ \big| \ \forall k \leq K : \, (l_k>0 \Rightarrow i^k_1 = b_k) \big\} \ ,
		\end{align*}
		then the above formula can be written as
		\begin{align}\label{EqGUE_Ex_Decomp1}
			& E_{Z}(l) = \sum\limits_{b_1,...,b_K \in [n]} \sum\limits_{\substack{\bm{i} \in \mathbb{I}_{l,b}}} \E\bigg[ \prod\limits_{k=1}^{K} Z_{i^k_1,i^k_2}  \cdots Z_{i^k_{l_k-1},i^k_{l_k}} \, Z_{i^k_{l_k},i^k_1} \bigg] \ .
		\end{align}
		The formula now also holds when some $l_i$ are zero, since then the summation over the corresponding $b_i$ brings the needed factors $n = \tr(\bm{Z}^0)$.\\
		\\
		Since complex standard normal random variables satisfy $\E[X^a \ol{X^b}] = \mathbbm{1}_{a=b} \, a!$ and real standard normal random variables satisfy $\E[X^a] = \mathbbm{1}_{a \text{ even}} \, (a-1)!!$ (with the convention $(-1)!!=1$), the mean in (\ref{EqGUE_Ex_Decomp1}) has the form
		\begin{align}
			& \E\bigg[ \prod\limits_{k=1}^{K} Z_{i^k_1,i^k_2}  \cdots Z_{i^k_{l_k-1},i^k_{l_k}} \, Z_{i^k_{l_k},i^k_1} \bigg] \nonumber\\
			& = \mathbbm{1}_{\substack{A(\bm{i}) \text{ is symmetric} \\ \text{\& diag. entries even}}} \bigg(\prod\limits_{\substack{v,w \in [n] \\ v < w}} A_{v,w}(\bm{i})!\bigg) \bigg(\prod\limits_{\substack{v \in [n]}} (A_{v,v}(\bm{i})-1)!!\bigg) =: \cW(A(\bm{i})) \ , \label{EqGUE_DefWeight}
		\end{align}
		where $A(\bm{i}) := \sum\limits_{k=1}^{K} A(G_{\bm{i}^k})$. Let $c := i^{1}_{l_1}$ - here the assumption $l_1>0$ is needed - and define the symmetric matrix
		\begin{align}\label{EqGUE_DefM}
			& M_{b,c} := \big(\mathbbm{1}_{i=b,j=c} + \mathbbm{1}_{i=c,j=b}\big)_{i,j \in [n]} \ .
		\end{align}
		For $c \neq b_1$, the above definition of $\cW(A(\bm{i}))$ then yields
		\begin{align*}
			& \cW(A(\bm{i})) = \cW(A(\bm{i})) - M_{b_1,c}) \, A_{c,b_1}(\bm{i})\\
			& = \cW(A(\bm{i}) - M_{b_1,c}) \bigg( A_{c,b_1}(G_{\bm{i}^1}) + \sum\limits_{s=2}^{K} A_{c,b_1}(G_{\bm{i}^s}) \bigg) \ ,
		\end{align*}
		where $A_{c,b_1}(\bm{i})$ is at least $1$, since it counts the number of occurrences of the connection $(i^1_{l_1},i^1_1) = (c,b_1)$. The fact that $\cW(A)$ is zero as soon as $A$ is not symmetric, allows us to swap $b_1$ and $c$ in arbitrary places of the formula, so we may also write
		\begin{align*}
			& \cW(A(\bm{i})) = \cW(A(\bm{i}) - M_{b_1,c}) \bigg( A_{b_1,c}(G_{\bm{i}^1}) + \sum\limits_{s=2}^{K} A_{b_1,c}(G_{\bm{i}^s}) \bigg) \ .
		\end{align*}
		Similarly, when $c=b_1$, the weight satisfies
		\begin{align*}
			& \cW(A(\bm{i})) = \cW(A(\bm{i})) - M_{b_1,b_1}) \, (A_{b_1,b_1}(\bm{i})-1)\\
			& = \cW(A(\bm{i}) - M_{b_1,b_1}) \bigg( A_{b_1,b_1}(G_{\bm{i}^1}) - 1 + \sum\limits_{s=2}^{K} A_{b_1,b_1}(G_{\bm{i}^s}) \bigg)
		\end{align*}
		and we can use these properties of the weight to decompose equality (\ref{EqGUE_Ex_Decomp1}) into
		\begin{align}
			& E_{Z}(l) = \overbrace{\sum\limits_{\substack{b \in [n]^K \\c \in [n]}} \sum\limits_{\substack{\bm{i} \in \mathbb{I}_{l,b} \\ i^1_{l_1}=c}} \cW(A(\bm{i}) - M_{b_1,c}) (A_{b_1,c}(G_{\bm{i}^1})-\mathbbm{1}_{b_1=c})}^{=: E^1_{Z}(l)} \nonumber\\
			& \hspace{1cm} + \underbrace{\sum\limits_{\substack{b \in [n]^K \\ c \in [n]}} \sum\limits_{\substack{\bm{i} \in \mathbb{I}_{l,b} \\ i^1_{l_1}=c}} \cW(A(\bm{i}) - M_{b_1,c}) \sum\limits_{s=2}^{K} A_{b_1,c}(G_{\bm{i}^s})}_{E^2_Z(l)} \ . \label{EqGUE_Ex_Decomp2}
		\end{align}
		For $E^1_{Z}(l)$ observe that $A_{b_1,c}(G_{\bm{i}^1})$ counts the number of occurrences of $(b_1,c)$ in $\bm{i}^1 = (i^1_1,...,i^1_{l_1})$ and, if this number is $\mathbbm{1}_{b_1=c}$, such an $\bm{i}$ will not contribute to $E^1_{Z}(l)$. Define
		\begin{align}\label{EqGUE_DefQ1}
			& Q^1_{b,c}(\bm{i}^1) = \{ q < l_1 \mid i^1_q=b, \, i^1_{q+1}=c \} \ ,
		\end{align}
		then $A_{b_1,c}(G_{\bm{i}^1}) - \mathbbm{1}_{b_1=c} = \#Q^1_{b_1,c}(\bm{i}^1)$. For each $q \in Q^1_{b_1,c}(\bm{i}^1)$ we give a unique way of splitting up $\bm{i}^1$ into two new routes $\langle\bm{i}^1\rangle_{q,1}$ and $\langle\bm{i}^1\rangle_{q,2}$. Define
		\begin{align*}
			& \langle\bm{i}^1\rangle_{q,1} := (\underbrace{i^1_1}_{=b_1},...,i^1_{q-1}) \in [n]^{q-1} \ \ ; \ \ \langle\bm{i}^1\rangle_{q,2} := (i^1_{q+2},...,\underbrace{i^1_{l_1}}_{=c}) \in [n]^{l_1-q-1}
		\end{align*}
		and observe that when going from the multi-graph $G_{\bm{i}^1}$ to $G_{\langle\bm{i}^1\rangle_{q,1}} \cup G_{\langle\bm{i}^1\rangle_{q,1}}$ we have only lost one edge in each direction (two edges without direction in the case of $b_1=c$) between the vertices $b_1$ and $c$ while all other edges remain. For
		\begin{align}\label{EqGUE_DefJ1}
			& J^1_{q}(\bm{i}) := \big( \langle\bm{i}^1\rangle_{q,1},\langle\bm{i}^1\rangle_{q,2},\bm{i}^2,...,\bm{i}^{K} \big)
		\end{align}
		it follows that
		\begin{align*}
			& A(J^1_{q}(\bm{i})) = A(\bm{i}) - M_{b_1,c} \ .
		\end{align*}
		Since we sum over all choices of $c$ and all the other elements in $\bm{i}^1$, which were not specified can still be chosen freely, this for $E^1_{l}$ implies
		\begin{align}
			& E^1_{Z}(l) = \sum\limits_{\substack{b \in [n]^{K} \\ c \in [n]}} \sum\limits_{\substack{\bm{i} \in \mathbb{I}_{l,b} \\ i^1_{l_1}=c}} \sum\limits_{q \in Q^1_{b_1,c}(\bm{i}^1)} \cW(A(J^1_q(\bm{i}))) \nonumber\\
			& = \sum\limits_{q=1}^{l_1-1} \sum\limits_{\substack{b \in [n]^{K} \\ c \in [n]}} \sum\limits_{\substack{\bm{i} \in \mathbb{I}_{l,b} \\ i^1_{l_1}=c \\ q \in Q^1_{b_1,c}(\bm{i}^1)}} \cW(A(J^1_q(\bm{i}))) \nonumber\\
			& = \sum\limits_{q=1}^{l_1-1} \sum\limits_{\substack{\tb \in [n]^{K+1}}} \sum\limits_{\substack{\bm{i} \in \mathbb{I}_{S_{2,q}(l),\tb}}} \cW(A(\bm{i})) \ , \label{EqGUE_E1Decomp1}
		\end{align}
		where
		\begin{align*}
			& S_{2,q}(l) := \big(q-1,l_1-q-1,l_2,...,l_{K}\big)
		\end{align*}
		is defined in analogy to $J^1_q(\bm{i})$. We have thus shown
		\begin{align}\label{EqGUE_E1Decomp2}
			& E^1_{Z}(l) = \sum\limits_{q=1}^{l_1-1} E_{Z}(S_{2,q}(l)) \ .
		\end{align}
		Next, for $E^2_{Z}(l)$ observe that $\sum\limits_{k=2}^{K} A_{b_1,c}(G_{\bm{i}^k})$ counts the number of occurrences of $(b_1,c)$ in
		\begin{align*}
			& (i^2_1,...,i^2_{l_2}), \cdots ,(i^{K}_1,...,i^{K}_{l_{K}})
		\end{align*}
		and define the set
		\begin{align}\label{EqGUE_DefQ2}
			& Q^2_{b,c}(\bm{i}) := \big\{ (k,q) \in \{2,...,K\} \times \N \ \big| \ q \leq l_k \text{ and } \nonumber\\
			& \hspace{2cm} q < l_k \Rightarrow (i^k_q = b \land i^k_{q+1}=c) , \, q = l_k \Rightarrow (i^k_{l_k} = b \land i^k_{1}=c) \big\} \ .
		\end{align}
		Again $\sum\limits_{k=2}^{K} A_{b_1,c}(G_{\bm{i}^k}) = \#Q^2_{b_1,c}(\bm{i})$ and this time for each $(k,q) \in Q^2_{b_1,c}(\bm{i})$ we give a unique way of joining the routes $\bm{i}^1$ and $\bm{i}^k$ into a new route
		\begin{align*}
			& \langle \bm{i}^1,\bm{i}^k \rangle_q := \begin{cases}
				\big( \underbrace{i^1_1}_{=b_1},...,\underbrace{i^1_{l_1}}_{=c}, i^k_{q+2},...,i^{k}_{l_k}, \underbrace{i^k_1}_{=b_k},...,i^{k}_{q-1} \big) \in [n]^{l_1+l_s-1} & \text{, if } q<l_k\\
				\big( \underbrace{i^1_1}_{=b_1},...,\underbrace{i^1_{l_1}}_{=c}, i^k_2,...,i^{k}_{l_k-1} \big) \in [n]^{l_1+l_s-1} & \text{, if } q=l_k\\
			\end{cases} \ .
		\end{align*}
		Going from the multigraph $G_{\bm{i}^1} \cup G_{\bm{i}^k}$ to $G_{\langle \bm{i}^1,\bm{i}^k \rangle_q}$ we only loose one edge in each direction between the vertices $b_1$ and $c$ (two undirected edges in the case $b_1=c$) while all other edges remain. For
		\begin{align}\label{EqGUE_DefJ2}
			& J^2_{(k,q)}(\bm{i}) := \big( \langle \bm{i}^1,\bm{i}^k \rangle_q, \bm{i}^2...,\widehat{\bm{i}^{k}},...,\bm{i}^{K} \big)
		\end{align}
		it follows that
		\begin{align*}
			& A(J^2_{(k,q)}(\bm{i})) = A(\bm{i}) - M_{b_1,c}
		\end{align*}
		and similarly to (\ref{EqGUE_E1Decomp1}) we get
		\begin{align}
			& E^2_{Z}(l) = \sum\limits_{\substack{b \in [n]^{K} \\ c \in [n]}} \sum\limits_{\substack{\bm{i} \in \mathbb{I}_{l,b} \\ i^1_{l_1}=c}} \sum\limits_{(k,q)\in Q^2_{b_1,c}} A(J^2_{(k,q)}(\bm{i})) \nonumber\\
			& = \sum\limits_{k=2}^K \sum\limits_{q=1}^{l_k} \sum\limits_{\substack{b \in [n]^{K} \\ c \in [n]}} \sum\limits_{\substack{\bm{i} \in \mathbb{I}_{l,b} \\ i^1_{l_1}=c \\ (k,q)\in Q^2_{b_1,c}}} A(J^2_{(k,q)}(\bm{i})) \nonumber\\
			& = \sum\limits_{k=2}^K \sum\limits_{q=1}^{l_k} \sum\limits_{\substack{\tb \in [p]^{K-1}}} \sum\limits_{\substack{\bm{i} \in \mathbb{I}_{S_{3,k}(l),\tb}}} A(\bm{i}) \nonumber\\
			& = \sum\limits_{k=2}^{K} \sum\limits_{r=1}^{l_k} E_{Z}(S_{3,k}(l)) = \sum\limits_{k=2}^{K} l_k E_{Z}(S_{3,k}(l)) \ , \label{EqGUE_E2Decomp}
		\end{align}
		where
		\begin{align*}
			& S_{3,k}(l) := \big(l_1+l_k-2,l_2,...,\widehat{l_k},...,l_{K}\big)
		\end{align*}
		is defined in analogy to $J^2_{(k,q)}$.\\
		\\
		Equalities (\ref{EqGUE_Ex_Decomp2}), (\ref{EqGUE_E1Decomp2}) and (\ref{EqGUE_E2Decomp}) together show
		\begin{align*}
			& E_Z(l) = \sum\limits_{q=1}^{l_1-1} E_{Z}(S_{2,q}(l)) + \sum\limits_{k=2}^{K} l_k E_{Z}(S_{3,k}(l)) \ . \qedhere
		\end{align*}
	\end{proof}

	\subsection{Remark (Translation to polygon gluings)}
	In terms of $\varepsilon_g(l)$ from (\ref{Eq_GUEOrderExpansion}) the recursion can be written as
	\begin{align*}
		& \sum\limits_{g=0}^{\lfloor g_{L,K}(0) \rfloor} n^{V_{L,K}(g)} \varepsilon_{g}(l)\\
		& = \sum\limits_{q=1}^{l_1-1} \sum\limits_{g=0}^{\lfloor g_{L-2,K+1}(0) \rfloor} n^{V_{L-2,K+1}(g)} \varepsilon_{g}(S_{2,q}(l)) + \sum\limits_{k=2}^{K} l_k \sum\limits_{g=0}^{\lfloor g_{L-2,K-1}(0) \rfloor} n^{V_{L-2,K-1}(g)} \varepsilon_{g}(S_{3,k}(l)) \ .
	\end{align*}
	With
	\begin{align*}
		& V_{L,K}(g) = 2-2g+\frac{L}{2}-K = V_{L-2,K+1}(g) = V_{L-2,K-1}(g-1)
	\end{align*}
	it implies
	\begin{align*}
		& \varepsilon_{g}(l) = \sum\limits_{q=1}^{l_1-1} \varepsilon_{g}(S_{2,q}(l)) + \sum\limits_{k=2}^{K} l_k \varepsilon_{g-1}(S_{3,k}(l)) \ ,
	\end{align*}
	when $l_1>0$, and $\varepsilon_{g}(0,l_2,...,l_K) = \varepsilon_{g}(l_2,...,l_K)$ otherwise.

	\section{Wishart Ensemble}\label{SectionWishart}
	For ease of notation in the following proofs we write the matrices $\bm{X}$ and $\bm{Y}$ as
	\begin{align*}
		& \bm{X} =: (X_{i,j})_{i \in [p], j \in [p+n] \setminus [p]} \ \text{ and } \ \bm{Y} =: (Y_{i,j})_{i \in [p], j \in [p+n] \setminus [p]} \ .
	\end{align*}
	
	\subsection{Theorem (Recursion for mixed trace moments of complex isotropic Wishart matrices)}\label{Thm_WC_Recursion}
	For $K \in \N$ and any $l=(l_1,...,l_{K}) \in \N_0^{K}$ the following recursive properties hold.
	\begin{itemize}
		\item[a)]
		If $l_1=0$, then $E_{Y}(l) = pE_{Y}(\cS_0(l))$ for
		\begin{align}\label{EqWC_DefS0}
			& \cS_0(0,l_2,...,l_K) := (l_2,...,l_K) \ .
		\end{align}
		
		\item[b)]
		If $l_1>0$, then
		\begin{align*}
			& E_Y(l) = nE_{Y}(\cS_1(l)) + \sum\limits_{r=1}^{l_1-1} E_{Y}(\cS_{2,r}(l)) + \sum\limits_{k=2}^{K} l_k E_{Y}(\cS_{3,l_k}(l)) \ ,
		\end{align*}
		where
		\begin{align}
			& \cS_{1}(l) := \big(l_1-1,l_2,...,l_{K}\big) \label{EqWC_DefS1}\\
			& \cS_{2,r}(l) := \big(r,l_1-r-1,l_2,...,l_{K}\big) \label{EqWC_DefS2}\\
			& \cS_{3,k,q}(l) := \big(l_1+l_k-1,l_2,...,\widehat{l_k},...,l_{K}\big) \ . \label{EqWC_DefS3}
		\end{align}
	\end{itemize}
	Since (a) reduces $K$ by one and (b) reduces $L := \sum\limits_{k=1}^K l_k$ by one, the end point
	\begin{align*}
		& E_{Y}\big( ((,)) \big) = \E\big[\tr\big((\bm{Y}\bm{Y}^*)^0\big)\big] = p
	\end{align*}
	for $(K,L) = (1,0)$ allows for the calculation of all mixed moments.
	\begin{proof}\
		\\
		The property (a) holds trivially, since the empty matrix product is here the $(p \times p)$ identity matrix. The proof of property (b) will take most of the section.\\
		\\
		For now assume $l_1,...,l_{K}>0$, then by expanding the sums in $\tr(\bm{S}^{l_i})$ we see
		\begin{align*}
			& E_{Y}(l) = \sum\limits_{\substack{\bm{i}^1 \in \mathcal{B}_{p,n,l_1} \\ : \\ \bm{i}^{K} \in \mathcal{B}_{p,n,l_{K}}}} \E\bigg[ \prod\limits_{k=1}^{K} \big( Y_{i^k_1,i^k_2} \ol{Y_{i^k_3,i^k_2}} \big) \cdots \big(Y_{i^k_{2l_k-1},i^k_{2l_k}} \ol{Y_{i^k_1,i^k_{2l_k}}} \big) \bigg] \ ,
		\end{align*}
		where
		\begin{align*}
			& \mathcal{B}_{p,n,l} := \big\{ (i_1,...,i_{2l}) \in \N^{2l} \ \big| \ \{i_1,i_3,...,i_{2l-1}\} \in [p] , \, \{i_2,i_4,...,i_{2l}\} \in [p+n]\setminus[p] \big\} \ .
		\end{align*}
		Extracting the first entries of each $\bm{i}^k$ and summing over them first yields
		\begin{align*}
			& E_{Y}(l) = \sum\limits_{b_1,...,b_K \in [p]} \sum\limits_{\substack{\bm{i}^1 \in \mathcal{B}_{p,n,l_1} \\ : \\ \bm{i}^{K} \in \mathcal{B}_{p,n,l_{K}} \\ \forall k \leq K: \, i^k_1 = b_k}} \E\bigg[ \prod\limits_{k=1}^{K} \big( Y_{i^k_1,i^k_2} \ol{Y_{i^k_3,i^k_2}} \big) \cdots \big( Y_{i^k_{2l_k-1},i^k_{2l_k}} \ol{Y_{i^k_1,i^k_{2l_k}}} \big) \bigg] \ .
		\end{align*}
		For ease of notation define
		\begin{align*}
			& \mathbb{B}_{l,b} := \big\{ \bm{i} = (\bm{i}^1,...,\bm{i}^K) \in \bigtimes\limits_{k=1}^K \mathcal{B}_{p,n,l_k} \ \big| \ \forall k \leq K: \, (l_k>0 \Rightarrow i^k_1 = b_k) \big\} \ ,
		\end{align*}
		then the above formula for $E_Y(l)$ becomes
		\begin{align}\label{EqWC_El_Decomp1}
			& E_{Y}(l) = \sum\limits_{b_1,...,b_K \in [p]} \sum\limits_{\bm{i} \in \mathbb{B}_{l,b}} \E\bigg[ \prod\limits_{k=1}^{K} \big( Y_{i^k_1,i^k_2} \ol{Y_{i^k_3,i^k_2}} \big) \cdots \big( Y_{i^k_{2l_k-1},i^k_{2l_k}} \ol{Y_{i^k_1,i^k_{2l_k}}} \big) \bigg] \ .
		\end{align}
		Just as (\ref{EqGUE_Ex_Decomp1}) in the GUE case, this formula now also holds for any $l_1,...,l_K \in \N_0^{K}$. This time the mean is seen to have the form
		\begin{align}
			& \E\bigg[ \prod\limits_{k=1}^{K} Y_{i^k_1,i^k_2} \ol{Y_{i^k_3,i^k_2}} \cdots Y_{i^k_{2l_k-1},i^k_{2l_k}} \ol{Y_{i^k_1,i^k_{2l_k}}} \bigg] \nonumber\\
			& = \mathbbm{1}_{A(\bm{i}) \text{ is symmetric}} \prod\limits_{\substack{v \in [p] \\ w \in [p+n] \setminus [p]}} A_{v,w}(\bm{i})! =: \cW(A(\bm{i})) \ , \label{EqWC_DefWeight}
		\end{align}
		where again $A(\bm{i}) := \sum\limits_{k=1}^{K} A(G_{\bm{i}^k})$. Let $c := i^{1}_{2l_1}$ and define the symmetric matrix
		\begin{align}\label{EqWC_DefM}
			& M_{b,c} := (\mathbbm{1}_{i=b,j=c \text{ or } i=c,j=b})_{i,j \in [p+n]} \ .
		\end{align}
		With the very same argument as the one leading up to (\ref{EqGUE_Ex_Decomp2}) - except that this time only the case $b_1 \neq c$ is needed - we can decompose equality (\ref{EqWC_El_Decomp1}) into
		\begin{align}
			& E_{Y}(l) = \overbrace{\sum\limits_{\substack{b \in [p]^K \\c \in [p+n] \setminus [p]}} \sum\limits_{\substack{\bm{i} \in \mathbb{B}_{l,b} \\ i^1_{2l_1}=c}} \cW(A(\bm{i}) - M_{b_1,c}) A_{b_1,c}(G_{\bm{i}^1})}^{=: E^1_{Y}(l)} \nonumber\\
			& \hspace{1cm} + \underbrace{\sum\limits_{\substack{b \in [p]^K \\c \in [p+n] \setminus [p]}} \sum\limits_{\substack{\bm{i} \in \mathbb{B}_{l,b} \\ i^1_{2l_1}=c}} \cW(A(\bm{i}) - M_{b_1,c}) \sum\limits_{s=k}^{K} A_{b_1,c}(G_{\bm{i}^k})}_{=: E^2_Y(l)} \ . \label{EqWC_Ex_Decomp2}
		\end{align}
		For $E^1_{Y}(l)$ observe that $A_{b_1,c}(G_{\bm{i}^1})$ counts the number of occurrences of $(b_1,c)$ in $\bm{i}^1 = (i^1_1,...,i^1_{2l_1})$ and, if this number is zero, such an $\bm{i}$ will not contribute to $E^1_{l}$. Define
		\begin{align}\label{EqWC_DefQ1}
			& Q^1_{b,c}(\bm{i}^1) = \{q<2l_1 \mid i^1_q=b, \, i^1_{q+1}=c \} \ ,
		\end{align}
		then clearly $A_{b_1,c}(G_{\bm{i}^1}) = \#Q^1_{b_1,c}(\bm{i}^1)$ and by the bipartite definition of $\mathcal{B}_{p,n,l_1}$ all elements in $Q^1_{b,c}(\bm{i}^1)$ must be odd. For each $q \in Q^1_{b_1,c}(\bm{i}^1)$ we give a unique way of splitting up $\bm{i}^1$ into two new routes $\langle\bm{i}^1\rangle_{q,1}$ and $\langle\bm{i}^1\rangle_{q,2}$. Define
		\begin{align*}
			& \langle\bm{i}^1\rangle_{q,1} := (i^1_{q+2},...,\underbrace{i^1_{2l_1}}_{=c}) \in \mathcal{B}_{p,n,l_1-\frac{q+1}{2}} \ \ ; \ \ \langle\bm{i}^1\rangle_{q,2} := (\underbrace{i^1_1}_{=b_1},...,i^1_{q-1}) \in \mathcal{B}_{p,n,\frac{q-1}{2}}
		\end{align*}
		and observe that between the multi-graphs $G_{\bm{i}^1}$ and $G_{\langle\bm{i}^1\rangle_{q,1}} \cup G_{\langle\bm{i}^1\rangle_{q,1}}$ we have only lost one edge in each direction between the vertices $b_1$ and $c$ while all other edges remain. For
		\begin{align}\label{EqWC_DefJ1}
			& J^1_{q}(\bm{i}) := \big( \langle\bm{i}^1\rangle_{q,1},\langle\bm{i}^1\rangle_{q,2},\bm{i}^2,...,\bm{i}^{K} \big)
		\end{align}
		it follows that
		\begin{align*}
			& A(J^1_{q}(\bm{i})) = A(\bm{i}) - M_{b_1,c} \ .
		\end{align*}
		Since we sum over all choices of $c$ and all the other elements in $\bm{i}^1$, which were not specified can still be chosen freely according the rules of $\mathbb{B}_{l,b}$, this for $E^1_{l}$ implies
		\begin{align}
			& E^1_{Y}(l) = \sum\limits_{\substack{b \in [p]^{K} \\ c \in [p+n] \setminus [p]}} \sum\limits_{\substack{\bm{i} \in \mathbb{B}_{l,b} \\ i^1_{2l_1}=c}} \sum\limits_{q \in Q^1_{b_1,c}(\bm{i}^1)} \cW(A(J^1_q(\bm{i}))) \nonumber\\
			& = \sum\limits_{\substack{q=1 \\ \text{odd}}}^{2l_1-1} \sum\limits_{\substack{b \in [p]^{K} \\ c \in [p+n] \setminus [p]}} \sum\limits_{\substack{\bm{i} \in \mathbb{B}_{l,b} \\ i^1_{2l_1}=c \\ q \in Q^1_{b_1,c}(\bm{i}^1)}} \cW(A(J^1_q(\bm{i}))) \nonumber\\
			& \overset{r=l_1-\frac{q-1}{2}}{=} \sum\limits_{\substack{r=0}}^{l_1-1} \sum\limits_{\substack{b \in [p]^{K} \\ b' \in [p] \text{, if } r>0 \\ c \in [p+n] \setminus [p]}} \sum\limits_{\substack{\bm{j} \in \mathbb{B}_{S_{2,r}(l),(b',b_1,...,b_K)} \\ r > 0 \Rightarrow j^1_{2r}=c}} \cW(A(\bm{j})) \nonumber\\
			& = \frac{n}{p} \sum\limits_{\substack{\tb \in [p]^{K+1}}} \sum\limits_{\substack{\bm{j} \in \mathbb{B}_{S_{2,0}(l),\tb}}} \cW(A(\bm{j})) \nonumber\\
			& \hspace{1cm} + \sum\limits_{\substack{r=1}}^{l_1-1} \sum\limits_{\substack{\tb \in [p]^{K+1}}} \sum\limits_{\substack{\bm{j} \in \mathbb{B}_{S_{2,r}(l),\tb}}} \cW(A(\bm{j})) \nonumber\\
			& = n E_Y(\cS_1(l)) + \sum\limits_{\substack{r=1}}^{l_1-1} E_{Y}(\cS_{2,r}(l)) \ , \label{EqWC_E1Decomp}
		\end{align}
		where
		\begin{align*}
			& \cS_{2,r}(l) = \big(r,l_1-r-1,l_2,...,l_{K}\big)
		\end{align*}
		is defined in analogy to $J^1_q(\bm{i})$ and
		\begin{align*}
			& \cS_{1}(l) = \big(l_1-1,l_2,...,l_{K}\big) \ .
		\end{align*}
		The factor $n$ in the first summand comes from the summation over $c$, which no longer had an influence on the latter sum in the case $q=2l_1-1 \ \Leftrightarrow \ r=0$. The temporary factor $\frac{1}{p}$ in the first summand was added and removed with the summation over $b' = \tb_1$.\\
		\\
		Next for $E^2_{Y}(l)$ observe that $\sum\limits_{k=2}^{K_1} A_{b_1,c}(G_{\bm{i}^k})$ counts the number of occurrences of $(b_1,c)$ in
		\begin{align*}
			& (i^2_1,...,i^2_{2l_2}), \cdots ,(i^{K}_1,...,i^{K}_{2l_{K}})
		\end{align*}
		and define the set
		\begin{align}\label{EqWC_DefQ2}
			& Q^2_{b,c}(\bm{i}) := \{ (k,q) \in \{2,...,K\} \times \N \mid i^k_q = b, \, i^k_{q+1}=c \} \ .
		\end{align}
		Again $\sum\limits_{k=2}^{K} A_{b_1,c}(G_{\bm{i}^k}) = \#Q^2_{b_1,c}$ and all tuples $(k,q) \in Q^2_{b,c}(\bm{i})$ must have odd $q$. This time for each $(k,q) \in Q^2_{b_1,c}(\bm{i})$ we give a unique way of joining the routes $\bm{i}^1$ and $\bm{i}^k$ into a new route
		\begin{align*}
			& \langle \bm{i}^1,\bm{i}^k \rangle_q := \big( \underbrace{i^k_1}_{=b_k},...,i^{k}_{q-1}, \underbrace{i^1_1}_{=b_1},...,\underbrace{i^1_{2l_1}}_{=c}, i^k_{q+2},...,i^{k}_{2l_k} \big) \in \mathcal{B}_{p,n,l_1+l_s-1} \ .
		\end{align*}
		Again we between the multigraphs $G_{\bm{i}^1} \cup G_{\bm{i}^k}$ and $G_{\langle \bm{i}^1,\bm{i}^k \rangle_q}$ only loose one edge in each direction between the vertices $b_1$ and $c$ while all other edges remain. For
		\begin{align}\label{EqWC_DefJ2}
			& J^2_{(k,q)}(\bm{i}) := \big( \langle \bm{i}^1,\bm{i}^k \rangle_q,,\bm{i}^2...,\widehat{\bm{i}^k},...,\bm{i}^{K} \big)
		\end{align}
		it follows that
		\begin{align*}
			& A(J^2_{(k,q)}(\bm{i})) = A(\bm{i}) - M_{b_1,c}
		\end{align*}
		and similarly to (\ref{EqGUE_E2Decomp}) we get
		\begin{align}
			& E^2_{Y}(l) = \sum\limits_{\substack{b \in [p]^{K} \\ c \in [p+n] \setminus [p]}} \sum\limits_{\substack{\bm{i} \in \mathbb{B}_{l,b} \\ i^1_{2l_1}=c}} \sum\limits_{(k,q)\in Q^2_{b_1,c}} A(J^2_{(k,q)}(\bm{i})) \nonumber\\
			& \overset{r=\frac{q+1}{2}}{=} \sum\limits_{k=2}^{K} \sum\limits_{r=1}^{l_k} \sum\limits_{\substack{\tilde{b} \in [p]^{K-1}}} \sum\limits_{\substack{\bm{i} \in \mathbb{B}_{S_{3,k}(l),y,\tilde{b}}}} A(\bm{i},\bm{j}) = \sum\limits_{k=2}^{K} l_k E_{Y}(S_{3,k}(l)) \ , \label{EqWC_E2Decomp}
		\end{align}
		where
		\begin{align*}
			& \cS_{3,k}(l) := \big(l_1+l_k-1,l_2,...,\widehat{l_k},...,l_{K}\big)
		\end{align*}
		is defined in analogy to $J^2_{k,2r-1}$.\\
		\\
		Equalities (\ref{EqWC_Ex_Decomp2}), (\ref{EqWC_E1Decomp}) and (\ref{EqWC_E2Decomp}) together show
		\begin{align*}
			& E_Y(l) = nE_{Y}(\cS_1(l)) + \sum\limits_{r=1}^{l_1-1} E_{Y}(\cS_{2,r}(l)) + \sum\limits_{k=2}^{K} l_k E_{Y}(\cS_{3,l_k}(l)) \ . \qedhere
		\end{align*}
	\end{proof}\

	\subsection{Theorem (Recursion for mixed trace moments of real isotropic Wishart matrices)}\label{Thm_WR_Recursion}
	For $K \in \N$ and any $l=(l_1,...,l_{K}) \in \N_0^{K}$ the following recursive properties hold.
	\begin{itemize}
		\item[a)]
		If $l_1=0$, then $E_{Y}(l) = pE_{Y}(\cS_0(l))$ for $\cS_0(0,l_2,...,l_K) := (l_2,...,l_K)$.
		
		\item[b)]
		If $l_1>0$, then
		\begin{align*}
			& E_X(l) = (n+l_1-1) E_X(\cS_1(l)) + \sum\limits_{r=1}^{l_1-1} E_X(\cS_{2,r}(l)) + 2\sum\limits_{k=2}^K l_k E_X(\cS_{3,k}(l))
		\end{align*}
		where $\cS_1,\cS_{2,q}$ and $\cS_{3,k}$ are as defined in Theorem \ref{Thm_WC_Recursion}.
	\end{itemize}
	Again the end point
	\begin{align*}
		& E_{X}\big( ((,)) \big) = \E\big[\tr\big((\bm{X}\bm{X}^T)^0\big)\big] = p
	\end{align*}
	for $(K,L) = (1,0)$ allows for the calculation of all mixed moments.
	\begin{proof}\
		\\
		In complete analogy to the proof of Theorem \ref{Thm_WC_Recursion} know (a) to hold and for the proof of (b) in the case $l_1>0$ have
		\begin{align}\label{EqWR_El_Decomp1}
			& E_{X}(l) = \sum\limits_{b_1,...,b_K \in [p]} \sum\limits_{\bm{i} \in \mathbb{B}_{l,b}} \E\bigg[ \prod\limits_{k=1}^{K} \big( X_{i^k_1,i^k_2} X_{i^k_3,i^k_2} \big) \cdots \big( X_{i^k_{2l_k-1},i^k_{2l_k}} X_{i^k_1,i^k_{2l_k}} \big) \bigg]
		\end{align}
		for any $l_1,...,l_K \in \N_0^{K}$. The weight is instead
		\begin{align}
			& \E\bigg[ \prod\limits_{k=1}^{K} X_{i^k_1,i^k_2} X_{i^k_3,i^k_2} \cdots X_{i^k_{2l_k-1},i^k_{2l_k}} X_{i^k_1,i^k_{2l_k}} \bigg] \nonumber\\
			& = \mathbbm{1}_{\tilde{A}(\bm{i}) \text{ has even entries}} \prod\limits_{\substack{v \in [p] \\ w \in [p+n] \setminus [p]}} \big( \tilde{A}_{v,w}(\bm{i}) - 1 \big)!! =: \cW(\tilde{A}(\bm{i})) \ , \label{EqWR_DefWeight}
		\end{align}
		where $\tilde{A}(\bm{i}) = \sum\limits_{k=1}^K A(\tilde{G}_{\bm{i}})$. Let $c := i^{1}_{2l_1}$ and define the matrix
		\begin{align}\label{EqWR_DefM}
			& \tilde{M}_{b,c} := 2(\mathbbm{1}_{i=b,j=c \text{ or } i=c,j=b})_{i,j \in [p+n]} \ .
		\end{align}
		By the definition of the weight $\cW(\tilde{A}(\bm{i}))$ it for all $\bm{i} \in \mathbb{B}_{l,b}$ with $\tilde{A}(\bm{i})_{b_1,c}\geq 2$ follows that
		\begin{align*}
			& \cW(\tilde{A}(\bm{i})) = \cW(\tilde{A}(\bm{i})-\tilde{M}_{b_1,c}) \big( \tilde{A}(\bm{i})_{b_1,c} - 1 \big)\\
			& = \cW(\tilde{A}(\bm{i})-\tilde{M}_{b_1,c}) \bigg( A_{b_1,c}(\tilde{G}_{\bm{i}^1}) - 1 + \sum\limits_{k=2}^K A_{b_1,c}(\tilde{G}_{\bm{i}^k}) \bigg) \ .
		\end{align*}
		Again this can be used to decompose equality (\ref{EqWR_El_Decomp1}) into
		\begin{align}
			& E_{X}(l) = \overbrace{\sum\limits_{\substack{b \in [p]^K \\c \in [p+n] \setminus [p]}} \sum\limits_{\substack{\bm{i} \in \mathbb{B}_{l,b} \\ i^1_{2l_1}=c}} \cW(\tilde{A}(\bm{i}) - \tilde{M}_{b_1,c}) \big( A_{b_1,c}(\tilde{G}_{\bm{i}^1}) - 1 \big)}^{=: E^1_{X}(l)} \nonumber\\
			& \hspace{1cm} + \underbrace{\sum\limits_{\substack{b \in [p]^K \\c \in [p+n] \setminus [p]}} \sum\limits_{\substack{\bm{i} \in \mathbb{B}_{l,b} \\ i^1_{2l_1}=c}} \cW(\tilde{A}(\bm{i}) - \tilde{M}_{b_1,c}) \sum\limits_{k=2}^{K} A_{b_1,c}(\tilde{G}_{\bm{i}^k})}_{=: E^2_X(l)} \ , \label{EqWR_Ex_Decomp2}
		\end{align}
		since all $\bm{i} \in \mathbb{B}_{l,b}$ contributing to the sum must satisfy $\tilde{A}(\bm{i})_{b_1,c}\geq 2$. Define
		\begin{align}
			& Q^0_{b,c}(\bm{i}^1) = \{q<2l_1 \mid i^1_q=c , \, i^1_{q+1}=b \} \label{EqWR_DefQ1}\\
			& Q^1_{b,c}(\bm{i}^1) = \{q<2l_1 \mid i^1_q=b , \, i^1_{q+1}=c \} = Q^0_{c,b}(\bm{i}) \ ,
		\end{align}
		then $A_{b_1,c}(\tilde{G}_{\bm{i}^1})-1 = \#Q^0_{b_1,c}(\bm{i}^1) + \#Q^1_{b_1,c}(\bm{i}^1)$ and by construction all elements of $Q^0_{b_1,c}(\bm{i}^1)$ must be even and all elements of $Q^1_{b_1,c}(\bm{i}^1)$ must be odd. For each $q \in Q^1_{b_1,c}(\bm{i}^1)$ let $J^1_q(\bm{i})$ be just as in the proof of Theorem \ref{Thm_WC_Recursion} and we also have $\tilde{A}(J^1_q(\bm{i})) = \tilde{A}(\bm{i}) - \tilde{M}_{b_1,c}$. For each $q \in Q^0_{b_1,c}(\bm{i}^1)$ let
		\begin{align*}
			& \langle \bm{i}^1 \rangle_{q} = \big( \underbrace{i^1_1}_{=b_1} ,...,i^1_{q-1},\underbrace{i^1_{2l_1}}_{=c} ,i^1_{2l_1-1},...,i^1_{q+3},i^1_{q+2} \big) \in \mathcal{B}_{p,n,l_1-1}
		\end{align*}
		and
		\begin{align}\label{EqWR_DefJ0}
			& J^0_q := \big( \langle \bm{i}^1 \rangle_{q}, \bm{i}^2,...,\bm{i}^K \big) \ ,
		\end{align}
		then it also holds that $\tilde{A}(J^0_q(\bm{i})) = \tilde{A}(\bm{i}) - \tilde{M}_{b_1,c}$ and we can similarly to (\ref{EqWC_E1Decomp}) see
		\begin{align*}
			& E^1_X(l) = \sum\limits_{\substack{b \in [p]^K \\c \in [p+n] \setminus [p]}} \sum\limits_{\substack{\bm{i} \in \mathbb{B}_{l,b} \\ i^1_{2l_1}=c}} \sum\limits_{q \in Q^0_{b_1,c}(\bm{i})} \cW(\tilde{A}(J^0_q(\bm{i})))\\
			& \hspace{1cm} + \sum\limits_{\substack{b \in [p]^K \\c \in [p+n] \setminus [p]}} \sum\limits_{\substack{\bm{i} \in \mathbb{B}_{l,b} \\ i^1_{2l_1}=c}} \sum\limits_{q \in Q^1_{b_1,c}(\bm{i})} \cW(\tilde{A}(J^1_q(\bm{i}))) \ .
		\end{align*}
		For the second summand the exact same steps as in (\ref{EqWC_E1Decomp}) show that it must be equal to
		\begin{align*}
			& n E_X(\cS_1(l)) + \sum\limits_{r=1}^{l-1} E_X(\cS_{2,r}(l)) \ .
		\end{align*}
		The first summand can likewise be calculated to be
		\begin{align*}
			& \sum\limits_{\substack{q=2 \\ \text{even}}}^{2l_1-2} \sum\limits_{\substack{b \in [p]^K \\c \in [p+n] \setminus [p]}} \sum\limits_{\substack{\bm{i} \in \mathbb{B}_{l,b} \\ i^1_{2l_1}=c \\ q \in Q^0_{b_1,c}(\bm{i})}} \cW(\tilde{A}(J^0_q(\bm{i})))\\
			& = (l_1-1) \sum\limits_{\substack{b \in [p]^K}} \sum\limits_{\substack{\bm{j} \in \mathbb{B}_{\cS_1(l),b}}} \cW(\tilde{A}(\bm{j})) = (l_1-1) E_X(\cS_1(l))
		\end{align*}
		and we have shown
		\begin{align}\label{EqWR_E1Calc}
			& E^1_X(l) = (n+l_1-1) E_X(\cS_1(l)) + \sum\limits_{r=1}^{l_1-1} E_X(\cS_{2,r}(l)) \ .
		\end{align}
		The same idea can be used to calculate $E_X^2(l)$. First define
		\begin{align}
			& Q^2_{b,c}(\bm{i}) = \{ (k,q) \in \{2,...,K\} \times \N \mid q < 2l_k, \, i^k_q=b, \, i^k_{q+1}=c \}\\
			& Q^3_{b,c}(\bm{i}) = \{ (k,q) \in \{2,...,K\} \times \N \mid q \leq 2l_k, \, i^k_q=c, \, i^k_{(q \mod 2l_k)+1}=b \} \ , \label{EqWR_DefQ3}
		\end{align}
		then again $\sum\limits_{k=2}^K A_{b_1,c}(\tilde{G}_{\bm{i}^k}) = \#Q^2_{b_1,c}(\bm{i}) + \#Q^3_{b_1,c}(\bm{i})$ and all tuples $(k,q)$ from $Q^2_{b,c}(\bm{i})$ must have odd $q$ while all tuples $(k,q)$ from $Q^3_{b,c}(\bm{i})$ must have even $q$. For $J^2_{(k,q)}(\bm{i})$ as in the proof of Theorem \ref{Thm_WC_Recursion}, it holds that $\tilde{A}(J^2_{(k,q)}(\bm{i})) = \tilde{A}(\bm{i}) - \tilde{M}_{b_1,c}$. Similarly for $(k,q) \in Q^3_{b_1,c}(\bm{i})$ define
		\begin{align*}
			& [ \bm{i}^1, \bm{i}^k ]_{q} := \begin{cases}
				\big( \overbrace{i^1_1}^{=b_1},...,\overbrace{i^1_{2l_1}}^{=c},i^k_q,i^k_{q-1},...,i^k_1,i^k_{2l_k},i^k_{2l_k-1},...,i^k_{q+2} \big) & \text{, if } q<2l_k\\
				\big( \underbrace{i^1_1}_{=b_1},...,\underbrace{i^1_{2l_1}}_{=c},i^k_{2l_k},i^k_{2l_k-1},...,i^k_2 \big) & \text{, if } q=2l_k
			\end{cases}
		\end{align*}
		and
		\begin{align}\label{EqWR_DefJ3}
			& J^3_{(k,q)} = \big( [\bm{i}^1,\bm{i}^k]_q,\bm{i}^2,...,\widehat{\bm{i}^k},...,\bm{i}^K \big)
		\end{align}
		then we also have $\tilde{A}(J^3_{(k,q)}(\bm{i})) = \tilde{A}(\bm{i}) - \tilde{M}_{b_1,c}$ and can further decompose $E^2_X(l)$ into
		\begin{align*}
			& E^2_X(l) = \sum\limits_{\substack{b \in [p]^K \\c \in [p+n] \setminus [p]}} \sum\limits_{\substack{\bm{i} \in \mathbb{B}_{l,b} \\ i^1_{2l_1}=c}} \sum\limits_{(k,q) \in Q^2_{b_1,c}(\bm{i})} \cW(\tilde{A}(J^2_{(k,q)}(\bm{i})))\\
			& \hspace{1cm} + \sum\limits_{\substack{b \in [p]^K \\c \in [p+n] \setminus [p]}} \sum\limits_{\substack{\bm{i} \in \mathbb{B}_{l,b} \\ i^1_{2l_1}=c}} \sum\limits_{(k,q) \in Q^3_{b_1,c}(\bm{i})} \cW(\tilde{A}(J^3_{(k,q)}(\bm{i}))) \ .
		\end{align*}
		The first summand is by the same steps as in (\ref{EqWC_E2Decomp}) equal to $\sum\limits_{k=2}^K l_k E_X(S_{3,k}(l))$. Note that $(k,2l_k) \in Q^3_{b_1,c}$ is only possible, if $b_k=b_1$. This will ensure that we don't gain a factor of $n$ from the summation over $b_k$ in the case $q=2l_k$ of the below calculation. The second summand is calculated to be
		\begin{align*}
			& \sum\limits_{k=2}^K \sum\limits_{\substack{q=2 \\ \text{even}}}^{2l_k} \sum\limits_{\substack{b \in [p]^K \\c \in [p+n] \setminus [p]}} \sum\limits_{\substack{\bm{i} \in \mathbb{B}_{l,b} \\ i^1_{2l_1}=c \\ (k,q) \in Q^3_{b_1,c}(\bm{i})}} \cW(\tilde{A}(J^3_{(k,q)}(\bm{i})))\\
			& = \sum\limits_{k=2}^K \sum\limits_{\substack{q=2 \\ \text{even}}}^{2l_k} \sum\limits_{\substack{\tb \in [p]^{K-1}}} \sum\limits_{\substack{\bm{j} \in \mathbb{B}_{S_{3,k}(l),\tb}}} \cW(\tilde{A}(\bm{j}))\\
			& = \sum\limits_{k=2}^K l_K E_X(\cS_{3,k}(l))
		\end{align*}
		and we have shown
		\begin{align}\label{EqWR_E2Calc}
			& E^2_X(l) = 2\sum\limits_{k=2}^K l_k E_X(\cS_{3,k}(l)) \ .
		\end{align}
		The equalities (\ref{EqWR_Ex_Decomp2}), (\ref{EqWR_E1Calc}) and (\ref{EqWR_E2Calc}) together prove
		\begin{align*}
			& E_X(l) = (n+l_1-1) E_X(\cS_1(l)) + \sum\limits_{r=1}^{l_1-1} E_X(\cS_{2,r}(l)) + 2\sum\limits_{k=2}^K l_k E_X(\cS_{3,k}(l)) \ . \qedhere
		\end{align*}
	\end{proof}

	\section{Gaussian Orthogonal Ensemble}
	\subsection{Theorem (Recursion for mixed trace moments of the GOE)}\label{Thm_GOE_Recursion}
	For $K \in \N$ and any $l=(l_1,...,l_{K}) \in \N_0^{K}$ - with even $L := l_1+...+l_K$ - the following recursive properties hold.
	\begin{itemize}
		\item[a)]
		If $l_1=0$, then $E_{\tZ}(l) = nE_{\tZ}(S_0(l))$ for
		\begin{align}\label{EqGOE_DefS0}
			& S_0(0,l_2,...,l_K) := (l_2,...,l_K) \ .
		\end{align}
		
		\item[b)]
		If $l_1>0$, then
		\begin{align*}
			& E_{\tZ}(l) = (l_1-1) E_{\tZ}(S_{1}(l)) + \sum\limits_{q=1}^{l_1-1} E_{\tZ}(S_{2,q}(l)) + 2 \sum\limits_{k=2}^K l_k E_{\tZ}(S_{3,k}(l)) \ ,
		\end{align*}
		where
		\begin{align}
			& S_{1}(x) := \big(l_1-2,l_2,...,l_{K}\big) \label{EqGOE_DefS1}
		\end{align}
		and $S_{2,q}$ as well as $S_{3,k}$ are as defined in Theorem \ref{Thm_GUE_Recursion}.
	\end{itemize}
	As before in Theorem \ref{Thm_GUE_Recursion} the starting point
	\begin{align*}
		& E_{\tZ}\big( ((,)) \big) = \E[\tr(\tilde{\bm{Z}})^0] = n
	\end{align*}
	for $(K,L) = (1,0)$ allows for the calculation of all mixed moments when $L$ is even and for odd $L$ the mixed trace moment must be zero.
	\begin{proof}\
		\\
		Like before, only property (b) needs to be shown. In complete analogy to (\ref{EqGUE_Ex_Decomp1}) we have
		\begin{align}\label{EqGOE_Ex_Decomp1}
			& E_{\tZ}(l) = \sum\limits_{b_1,...,b_K \in [n]} \sum\limits_{\substack{\bm{i} \in \mathbb{I}_{l,b}}} \E\bigg[ \prod\limits_{k=1}^{K} {\tZ}_{i^k_1,i^k_2}  \cdots {\tZ}_{i^k_{l_k-1},i^k_{l_k}} \, {\tZ}_{i^k_{l_k},i^k_1} \bigg] \ .
		\end{align}
		This time the mean has the form
		\begin{align}\label{EqGOE_DefWeight}
			& \E\bigg[ \prod\limits_{k=1}^{K} {\tZ}_{i^k_1,i^k_2} \cdots {\tZ}_{i^k_{l_k-1},i^k_{l_k}} \, {\tZ}_{i^k_{l_k},i^k_1} \bigg] \nonumber\\
			& = \mathbbm{1}_{\substack{\tA(\bm{i}) \text{ has even entries}}} \bigg( \prod\limits_{\substack{v,w \in [n] \\ v<w}} (\tA_{v,w}(\bm{i})-1)!! \bigg) \bigg( \prod\limits_{v \in [n]} 2^{\frac{\tA_{v,v}(\bm{i})}{2}} (\tA_{v,v}(\bm{i})-1)!! \bigg) =: \cW(\tA(\bm{i})) \ ,
		\end{align}
		which implies that it can with
		\begin{align}\label{EqGOE_DefM}
			& \tilde{M}_{b,c} := 2(\mathbbm{1}_{i=b,j=c \text{ or } i=c, j=b})_{i,j \in [n]}
		\end{align}
		be decomposed into
		\begin{align*}
			& \cW(\tA(\bm{i})) = \cW(\tA(\bm{i}) - \tilde{M}_{b_1,c}) \, 2^{\mathbbm{1}_{b_1=c}} \big( \tA_{b_1,c}(\bm{i}) - 1 \big)\\
			& = \cW(\tA(\bm{i}) - \tilde{M}_{b_1,c}) \, 2^{\mathbbm{1}_{b_1=c}} \bigg( A_{b_1,c}(\tilde{G}_{\bm{i}^1}) - 1 + \sum\limits_{k=2}^K A_{b_1,c}(\tilde{G}_{\bm{i}^k}) \bigg) \ .
		\end{align*}
		As is by now standard, we thus change (\ref{EqGOE_Ex_Decomp1}) into
		\begin{align}\label{EqGOE_Ex_Decomp2}
			& E_{\tZ}(l) = \overbrace{\sum\limits_{\substack{b_1,...,b_K \in [n] \\ c \in [n]}} 2^{\mathbbm{1}_{b_1=c}} \sum\limits_{\substack{\bm{i} \in \mathbb{I}_{l,b} \\ i^1_{l_1}=c}} \cW(\tA(\bm{i}) - \tilde{M}_{b_1,c}) (A_{b_1,c}(\tilde{G}_{\bm{i}^1}) - 1)}^{E^1_{\tZ}(l)} \nonumber\\
			& \hspace{1cm} + \underbrace{\sum\limits_{\substack{b_1,...,b_K \in [n] \\ c \in [n]}} 2^{\mathbbm{1}_{b_1=c}} \sum\limits_{\substack{\bm{i} \in \mathbb{I}_{l,b} \\ i^1_{l_1}=c}} \cW(\tA(\bm{i}) - \tilde{M}_{b_1,c}) \sum\limits_{k=2}^K A_{b_1,c}(\tilde{G}_{\bm{i}^k})}_{=: E^2_{\tZ}(l)} \ .
		\end{align}
		For
		\begin{align}
			& Q^0_{b,c}(\bm{i}^1) := \{q < l_1 \mid i^1_q = c, \, i^1_{q+1}=b\} \label{EqGOE_DefQ1}\\
			& Q^1_{b,c}(\bm{i}^1) := \{q < l_1 \mid i^1_q = b, \, i^1_{q+1}=c\} = Q^0_{c,b}(\bm{i}^1)
		\end{align}
		it holds that $2^{\mathbbm{1}_{b_1=c}} (A_{b_1,c}(\tilde{G}_{\bm{i}^1}) - 1) = \#Q^0_{b_1,c}(\bm{i}^1) + \#Q^1_{b_1,c}(\bm{i}^1)$ and with $J^1_q(\bm{i})$ as in (\ref{EqGUE_DefJ1}) we for each $q \in Q^1_{b_1,c}(\bm{i}^1)$ have
		\begin{align*}
			& \tilde{A}(J^1_q(\bm{i})) = \tilde{A}(\bm{i}) - \tilde{M}_{b_1,c} \ .
		\end{align*}
		For each $q \in Q^0_{b_1,c}(\bm{i}^1)$ we similarly to (\ref{EqWR_DefJ0}) define
		\begin{align*}
			& \langle \bm{i}^1 \rangle_q = \big( \underbrace{i^1_1}_{=b_1^1},...,i^1_{q-1},\underbrace{i^1_{l_1}}_{=c},i^1_{l_1-1},...,i^1_{q+3},i^1_{q+2} \big) \in [n]^{l_1-2}
		\end{align*}
		and
		\begin{align}\label{EqGOE_DefJ0}
			& J^0_q := \big( \langle \bm{i}^1 \rangle_q,\bm{i}^2,...,\bm{i}^K \big) \ ,
		\end{align}
		then it also follows that $\tilde{A}(J^0_q(\bm{i})) = \tilde{A}(\bm{i}) - \tilde{M}_{b_1,c}$, since only the connections between $(i^1_{q},i^1_{q+1}) = (c,b)$ and $(i^1_{l_1},i^1_1) = (c,b)$ are removed by $J^0_q$. Using $J^0$ and $J^1$ we write $E^1_{\tZ}(l)$ as
		\begin{align*}
			& E^1_{\tZ}(l) = \sum\limits_{\substack{b_1,...,b_K \in [n] \\ c \in [n]}} \sum\limits_{\substack{\bm{i} \in \mathbb{I}_{l,b} \\ i^1_{l_1}=c}} \sum\limits_{q \in Q^0_{b_1,c}(\bm{i}^1)} \cW(\tA(J^0_q))\\
			& \hspace{1cm} + \sum\limits_{\substack{b_1,...,b_K \in [n] \\ c \in [n]}} \sum\limits_{\substack{\bm{i} \in \mathbb{I}_{l,b} \\ i^1_{l_1}=c}} \sum\limits_{q \in Q^1_{b_1,c}(\bm{i}^1)} \cW(\tA(J^1_q)) \ .
		\end{align*}
		For the second summand the exact same steps as in (\ref{EqGUE_E1Decomp1}) and (\ref{EqGUE_E1Decomp2}) show that it must be equal to
		\begin{align*}
			& \sum\limits_{q=1}^{l_1-1} E_{\tZ}(S_{2,q}(l)) \ .
		\end{align*}
		Combining the ideas of (\ref{EqGUE_E1Decomp1}) and (\ref{EqWR_E1Calc}) we for the first summand observe
		\begin{align*}
			& \sum\limits_{\substack{b_1,...,b_K \in [n] \\ c \in [n]}} \sum\limits_{\substack{\bm{i} \in \mathbb{I}_{l,b} \\ i^1_{l_1}=c}} \sum\limits_{q \in Q^0_{b_1,c}(\bm{i}^1)} \cW(\tA(J^0_q)) = \sum\limits_{q=1}^{l_1-1} \sum\limits_{\substack{b_1,...,b_K \in [n] \\ c \in [n]}} \sum\limits_{\substack{\bm{i} \in \mathbb{I}_{l,b} \\ i^1_{l_1}=c \\ q \in Q^0_{b_1,c}(\bm{i}^1)}} \cW(\tA(J^0_q))\\
			& = (l_1-1) \sum\limits_{\substack{b_1,...,b_K \in [n]}} \sum\limits_{\substack{\bm{j} \in \mathbb{I}_{S_1(l),b}}} \cW(\tA(\bm{j})) = (l_1-1) E_{\tZ}(S_{1}(l)) \ ,
		\end{align*}
		where
		\begin{align*}
			& S_1(l) := \big( l_1-2,l_2,...,l_K \big) \ .
		\end{align*}
		We have thus shown
		\begin{align}\label{EqGOE_E1Calc}
			& E^1_{\tZ}(l) = (l_1-1) E_{\tZ}(S_{1}(l)) + \sum\limits_{q=1}^{l_1-1} E_{\tZ}(S_{2,q}(l)) \ .
		\end{align}
		With analogous steps to those used for the proof of (\ref{EqWR_E2Calc}) we will now calculate $E^2_{\tZ}(l)$. First define
		\begin{align}
			& Q^2_{b,c}(\bm{i}) := \{(k,q) \in \{2,...,K\} \times \N \mid q \leq l_k, i^k_q=b, i^k_{(q \mod l_k)+1}=c\}\\
			& Q^3_{b,c}(\bm{i}) := \{(k,q) \in \{2,...,K\} \times \N \mid q \leq l_k, i^k_q=c, i^k_{(q \mod l_k)+1}=b\} = Q^2_{c,b}(\bm{i}) \ , \label{EqGOE_DefQ3}
		\end{align}
		then we have $2^{\mathbbm{1}_{b_1=c}} \sum\limits_{k=2}^K A_{b_1,c}(\tilde{G}_{\bm{i}^k}) = \#Q^2_{b_1,c}(\bm{i}) + \#Q^3_{b_1,c}(\bm{i})$. For $J^2_{(k,q)}(\bm{i})$ as in (\ref{EqGUE_DefJ2}) it holds that $\tilde{A}(J^2_{(k,q)}(\bm{i})) = \tilde{A}(\bm{i})-\tilde{M}_{b_1,c_1}$ and for each $(k,q) \in Q^3_{b_1,c}(\bm{i})$ we similarly to (\ref{EqWR_DefJ3}) define
		\begin{align*}
			& [ \bm{i}^1, \bm{i}^k ]_{q} := \begin{cases}
				\big( \overbrace{i^1_1}^{=b_1},...,\overbrace{i^1_{l_1}}^{=c},i^k_q,i^k_{q-1},...,i^k_1,i^k_{l_k},i^k_{l_k-1},...,i^k_{q+2} \big) & \text{, if } q<l_k\\
				\big( \underbrace{i^1_1}_{=b_1},...,\underbrace{i^1_{l_1}}_{=c},i^k_{l_k},i^k_{l_k-1},...,i^k_2 \big) & \text{, if } q=l_k
			\end{cases}
		\end{align*}
		and
		\begin{align}\label{EqGOE_DefJ3}
			& J^3_{(k,q)} = \big( [\bm{i}^1,\bm{i}^k]_q,\bm{i}^2,...,\widehat{\bm{i}^k},...,\bm{i}^K \big) \ .
		\end{align}
		Again by construction of $J^3$ it holds that $\tilde{A}(J^3_{(k,q)}(\bm{i})) = \tilde{A}(\bm{i})-\tilde{M}_{b_1,c_1}$ and we can decompose $E^2_{\tZ}(l)$ into
		\begin{align*}
			& E^2_{\tZ}(l) = \sum\limits_{\substack{b_1,...,b_K \in [n] \\ c \in [n]}} \sum\limits_{\substack{\bm{i} \in \mathbb{I}_{l,b} \\ i^1_{l_1}=c}} \sum\limits_{(k,q) \in Q^2_{b_1,c}(\bm{i})} \cW(\tilde{A}(J^2_{(k,q)}(\bm{i})))\\
			& \hspace{1cm} + \sum\limits_{\substack{b_1,...,b_K \in [n] \\ c \in [n]}} \sum\limits_{\substack{\bm{i} \in \mathbb{I}_{l,b} \\ i^1_{l_1}=c}} \sum\limits_{(k,q) \in Q^3_{b_1,c}(\bm{i})} \cW(\tilde{A}(J^3_{(k,q)}(\bm{i}))) \ .
		\end{align*}
		The first summand is by the same steps as in (\ref{EqGUE_E2Decomp}) equal to
		\begin{align*}
			& \sum\limits_{k=2}^K l_k E_{\tZ}(S_{3,k}(l))
		\end{align*}
		and for the second summand the steps are again analogous to (\ref{EqGUE_E2Decomp}) and (\ref{EqWR_E2Calc}) and we find
		\begin{align*}
			& \sum\limits_{\substack{b_1,...,b_K \in [n] \\ c \in [n]}} \sum\limits_{\substack{\bm{i} \in \mathbb{I}_{l,b} \\ i^1_{l_1}=c}} \sum\limits_{(k,q) \in Q^3_{b_1,c}(\bm{i})} \cW(\tilde{A}(J^3_{(k,q)}(\bm{i})))\\
			& = \sum\limits_{k=2}^K \sum\limits_{q=1}^{l_k} \sum\limits_{\substack{b \in [n]^K \\ c \in [n]}} \sum\limits_{\substack{\bm{i} \in \mathbb{I}_{l,b} \\ i^1_{l_1}=c \\ (k,q) \in Q^3_{b_1,c}(\bm{i})}} \cW(\tilde{A}(J^3_{(k,q)}(\bm{i})))\\
			& = \sum\limits_{k=2}^K \sum\limits_{q=1}^{l_k} \sum\limits_{\substack{\tilde{b} \in [n]^{K-1}}} \sum\limits_{\substack{\bm{i} \in \mathbb{I}_{l,\tilde{b}}}} \cW(\tilde{A}(\bm{j}))\\
			& = \sum\limits_{k=2}^K l_k E_{\tZ}(S_{3,k}(l)) \ ,
		\end{align*}
		so we have shown
		\begin{align}\label{EqGOE_E2Calc}
			& E^2_{\tZ}(l) = 2 \sum\limits_{k=2}^K l_k E_{\tZ}(S_{3,k}(l)) \ .
		\end{align}
		Finally equalities (\ref{EqGOE_Ex_Decomp2}), (\ref{EqGOE_E1Calc}) and (\ref{EqGOE_E2Calc}) together yield
		\begin{align*}
			& E_{\tZ}(l) = (l_1-1) E_{\tZ}(S_{1}(l)) + \sum\limits_{q=1}^{l_1-1} E_{\tZ}(S_{2,q}(l)) + 2 \sum\limits_{k=2}^K l_k E_{\tZ}(S_{3,k}(l)) \ .
		\end{align*}
	\end{proof}

	\section*{List of symbols}\label{ListOfSymbols}
	\begin{itemize}
		
		\item[] $A(G)$ \tabto{1.5cm} Adjacency matrix to a graph $G$ (see Definition \ref{Def_Graphs})
		
		\item[] $A(\bm{i})$ \tabto{1.5cm} Sum of adjacency matrices $A(G_{\bm{i}^1}),...,A(G_{\bm{i}^K})$ (see Definition \ref{Def_Graphs})
		
		\item[] $\tilde{A}(\bm{i})$ \tabto{1.5cm} Sum of adjacency matrices $A(\tilde{G}_{\bm{i}^1}),...,A(\tilde{G}_{\bm{i}^K})$ (see Definition \ref{Def_Graphs})
		
		\item[] $e^k_q$ \tabto{1.5cm} an edge in a graph/polygon (see Subsection \ref{Subsection_Polygongluings} or Definition \ref{Def_Graphs})
		
		\item[] $E_X$ \tabto{1.5cm} Mixed trace moment of $\bm{X}\bm{X}^*$ (see Introduction)
		
		\item[] $E_Y$ \tabto{1.5cm} Mixed trace moment of $\bm{Y}\bm{Y}^*$ (see Introduction)
		
		\item[] $E_Z$ \tabto{1.5cm} Mixed trace moment of $\bm{Z}$ (see Introduction)
		
		\item[] $E_{\tZ}$ \tabto{1.5cm} Mixed trace moment of $\tilde{\bm{Z}}$ (see Introduction)
		
		\item[] $\varepsilon_g$ \tabto{1.5cm} number of different ways to glue edges of polygons (see Subsection \ref{Subsection_Polygongluings})
		
		\item[] $g$ \tabto{1.5cm} genus of a surface resulting from a polygon-gluing (see Subsection \ref{Subsection_Polygongluings})
		
		\item[] $G_{\bm{i}}$ \tabto{1.5cm} a directed multi-graph such that $\bm{i}$ is an Euler-tour (see Definition \ref{Def_Graphs})
		
		\item[] $\tilde{G}_{\bm{i}}$ \tabto{1.5cm} an undirected multi-graph such that $\bm{i}$ is an Euler-tour (see Definition \ref{Def_Graphs})
		
		\item[] $J^0_{(k,q)}$ \tabto{1.5cm} a locally defined operator on route sequences (see (\ref{EqWR_DefJ0}) or (\ref{EqGOE_DefJ0}))
		
		\item[] $J^1_{(k,q)}$ \tabto{1.5cm} a locally defined operator on route sequences (see (\ref{EqGUE_DefJ1}) or (\ref{EqWC_DefJ1}))
		
		\item[] $J^2_{(k,q)}$ \tabto{1.5cm} a locally defined operator on route sequences (see (\ref{EqGUE_DefJ2}) or (\ref{EqWC_DefJ2}))
		
		\item[] $J^3_{(k,q)}$ \tabto{1.5cm} a locally defined operator on route sequences (see (\ref{EqWR_DefJ3}) or (\ref{EqGOE_DefJ3}))
		
		\item[] $K$ \tabto{1.5cm} number of components to a layout $l=(l_1,...,l_K) \in \N_0^K$ (see Introduction)
		
		\item[] $l$ \tabto{1.5cm} a 'layout' $l=(l_1,...,l_K) \in \N_0^K$ (see Introduction)
		
		\item[] $L$ \tabto{1.5cm} the sum $\sum\limits_{k=1}^K l_k$ to a layout $l=(l_1,...,l_K) \in \N_0^K$ (see Introduction)
		
		\item[] $M_{b,c}$ \tabto{1.5cm} describes the symmetric matrix $\big(\mathbbm{1}_{i=b,j=c} + \mathbbm{1}_{i=c,j=b}\big)_{i,j}$ (see (\ref{EqGUE_DefM}) or (\ref{EqWC_DefM}))
		
		\item[] $\tilde{M}_{b,c}$ \tabto{1.5cm} describes the symmetric matrix $2(\mathbbm{1}_{i=b,j=c \text{ or } i=c,j=b})_{i,j}$ (see (\ref{EqGOE_DefM}) or (\ref{EqWR_DefM}))
		
		\item[] $n$ \tabto{1.5cm} dimension of a random matrix (see Introduction)
		
		\item[] $p$ \tabto{1.5cm} dimension of a random matrix (see Introduction)
		
		\item[] $P_m$ \tabto{1.5cm} a polygon with $m$ edges (see Subsection \ref{Subsection_Polygongluings})
		
		\item[] $\pi$ \tabto{1.5cm} a pairing of edges of polygons (see Subsection \ref{Subsection_Polygongluings})
		
		\item[] $Q^{0}_{b,c}$ \tabto{1.5cm} a locally defined index-set (see (\ref{EqWR_DefQ1}) or (\ref{EqGOE_DefQ1}))
		
		\item[] $Q^{1}_{b,c}$ \tabto{1.5cm} a locally defined index-set (see (\ref{EqGUE_DefQ1}) or (\ref{EqWC_DefQ1}))
		
		\item[] $Q^{2}_{b,c}$ \tabto{1.5cm} a locally defined index-set (see (\ref{EqGUE_DefQ2}) or (\ref{EqWC_DefQ2}))
		
		\item[] $Q^{3}_{b,c}$ \tabto{1.5cm} a locally defined index-set (see (\ref{EqWR_DefQ3}) or (\ref{EqGOE_DefQ3}))
		
		\item[] $S_\bullet$ \tabto{1.5cm} an operator on layouts (see (\ref{EqGUE_DefS0})-(\ref{EqGUE_DefS3}) and (\ref{EqGOE_DefS1}))
		
		\item[] $\cS_\bullet$ \tabto{1.5cm} an operator on layouts (see (\ref{EqWC_DefS0})-(\ref{EqWC_DefS3}))
		
		\item[] $\operatorname{tr}$ \tabto{1.5cm} trace of a square matrix
		
		\item[] $V$ \tabto{1.5cm} number of vertices in certain graphs (see Subsection \ref{Subsection_Polygongluings} or Definition \ref{Def_Graphs})
		
		\item[] $\cW(A)$ \tabto{1.5cm} a locally defined weight (see (\ref{EqGUE_DefWeight}), (\ref{EqWC_DefWeight}), (\ref{EqWR_DefWeight}) or (\ref{EqGOE_DefWeight}))
		
		\item[] $\bm{X}$ \tabto{1.5cm} a $(p \times n)$ random matrix with iid $\mathcal{N}(0,1)$ entries (see Introduction)
		
		\item[] $\bm{Y}$ \tabto{1.5cm} a $(p \times n)$ random matrix with iid $\mathcal{CN}(0,1)$ entries (see Introduction)
		
		\item[] $\bm{Z}$ \tabto{1.5cm} an $(n \times n)$ random matrix with such that $\frac{1}{\sqrt{n}} \bm{Z}$ is of the GUE\\
		\phantom{.}\hspace{2cm}(see Introduction)
		
		\item[] $\tilde{\bm{Z}}$ \tabto{1.5cm} an $(n \times n)$ random matrix with such that $\frac{1}{\sqrt{n}} \tilde{\bm{Z}}$ is of the GOE\\
		\phantom{.}\hspace{2cm}(see Introduction)
	\end{itemize}


\end{document}